\documentclass[10pt,twoside]{article}

\usepackage{amsfonts}
\usepackage{amsmath}
\usepackage{amsthm}
\usepackage{amssymb}
\usepackage{dsfont}

\usepackage[latin1]{inputenc}

\newcommand{\N}{\mathbb{N}}
\newcommand{\Z}{\mathbb{Z}}

\newcommand{\R}{\mathbb{R}}
\newcommand{\C}{\mathbb{C}}

\newtheorem{theorem}{Theorem}[section]{\bf}{\it }
\newtheorem{proposition}{Proposition}[section]{\bf}{\it }
\newtheorem{lemma}{Lemma}[section]{\bf}{\it }
\newtheorem{remark}{Remark}[section]{\bf}{\it }
\newtheorem{definition}{Definition}[section]{\bf}{\it }
\newtheorem{assumption}{Assumption}[section]{\bf}{\it }
\newtheorem{corrollary}{Corollary}[section]{\bf}{\it }

\numberwithin{equation}{section}

\def \p{\partial}
\def \ds{\displaystyle}
\def \e{\varepsilon}
\def \one{\mathds{1}}
\def \bone{\mathbf{1}}
\def \bzero{\mathbf{0}}
\def \rmi{{\rm i}}

\title{On the asymptotic stability of small nonlinear Dirac standing
waves in a resonant case}

\begin{document}
\author{Nabile Boussaid\\
CEREMADE, 
Universit\'e Paris-Dauphine,\\ 
Place du Mar\'echal De Lattre De Tassigny 
F-75775 Paris C\'edex 16 - France\\
boussaid@ceremade.dauphine.fr}


\maketitle

\begin{abstract}
We study the behavior of perturbations of small nonlinear Dirac
standing waves. We assume that the linear Dirac operator of
reference $H=D_m+V$ has only two double eigenvalues and that
degeneracies are due to a symmetry of $H$ (theorem of Kramers). In
this case, we can build a small $4$-dimensional manifold of
stationary solutions tangent to the first eigenspace of $H$.

Then we assume that a resonance condition holds and we build a
center manifold of real codimension $8$ around each stationary
solution. Inside this center manifold any $H^{s}$ perturbation of
stationary solutions, with $s>2$, stabilizes towards a standing
wave. We also build center-stable and center-unstable manifolds each
one of real codimension $4$. Inside each of these manifolds, we obtain
stabilization towards the center manifold in one direction of time,
while in the other, we have instability. Eventually, outside all
these manifolds, we have instability in the two directions of time.

For localized perturbations inside the center manifold, we obtain a nonlinear 
scattering result.

\end{abstract}

\section*{Introduction}
We study the asymptotic stability of stationary solutions of a
time-dependent nonlinear Dirac equation.

A localized stationary solution of a given time-dependent equation
represents a bound state of a particle. Like Ra\~nada \cite{Ranada},
we call it a \emph{particle-like solution} (PLS). Many works have
been devoted to the proof of the existence of such solutions for a
wide variety of equations. Although their stability is a crucial
problem (in particular in numerical computation or experiment), a
smaller attention has been deserved to this issue.

In this paper, we deal with the problem of stability of small PLS of
the following nonlinear Dirac equation:
\begin{equation}
\rmi\p_t\psi=(D_m+V)\psi+\nabla F(\psi) \tag{NLDE}\label{nld_intro}
\end{equation}
where $\nabla F$ is the gradient of $F:\C^4\mapsto \R$ for the
standard scalar product of $\R^8$. Here,~$D_m$ is the usual Dirac
operator, see Thaller~\cite{Thaller}, acting on $L^2(\R^3,\C^4)$
\begin{equation*}
D_m=\alpha\cdot\left(-\rmi\nabla\right)+m\beta=
-\rmi\sum_{k=1}^3\alpha_k\partial_k + m\beta
\end{equation*}
where~$m\in \R_+^*$,
~$\alpha=\left(\alpha_1,\alpha_2,\alpha_3\right)$, $\beta$
are~$\C^4$ hermitian matrices satisfying:
\begin{equation*}
\begin{cases}
\alpha_i\alpha_k+\alpha_k\alpha_i=2\delta_{ik}\bone_{\C^4},
&i,k\in\{1,2,3\},\\
\alpha_i\beta+\beta\alpha_i=\bzero_{\C^4},&i\in\{1,2,3\},\\
\beta^2=\bone_{\C^4}.
\end{cases}
\end{equation*}
Here, we choose
\begin{gather*}
\alpha_i=\left(
\begin{array}{cc}
0&\sigma_i\\
\sigma_i&0
\end{array} \right)\quad\mbox{ and }\quad\beta=\left(\begin{array}{cc}
I_{\C^2}&0\\
0&-I_{\C^2}
\end{array} \right)\\
\mbox{ where  } \sigma_1=\left(
\begin{array}{cc}
0&\; 1\\
1&\;0
\end{array} \right),\quad \quad\sigma_2=\left(
\begin{array}{cc}
0&-\rmi\\
\rmi&0
\end{array} \right)\quad\mbox{ and } \quad\sigma_3=\left(
\begin{array}{cc}
1&0\\
0&-1
\end{array} \right).
\end{gather*}
In~\eqref{nld_intro},~$V$ is the external potential field
and~$F:\C^4\mapsto\R$ is a nonlinearity with the following gauge
invariance:
\begin{equation}                                                \label{Eq:GaugeInvariance}
\forall(\theta,z)\in\R\times\C^4,\quad F(e^{i\theta}z)=F(z).
\end{equation}
Some additional assumptions on~$F$ and $V$ will be made in the
sequel. Stationary solutions (PLS) of \eqref{nld_intro} take the
form $\psi(t,x)=e^{-\rmi Et}\phi(x)$ where~$\phi$ satisfies
\begin{equation}
E\phi=(D_m+V)\phi+\nabla F(\phi). \tag{PLSE}\label{eq:stationary}
\end{equation}
We show that there exists a manifold of small solutions to
\eqref{eq:stationary} tangent to the first eigenspace of $D_m+V$
(see Proposition \ref{Prop:ManifoldPLS} below).

\medskip

In the Schrödinger case, orbital stability results (see e.g
\cite{CazenaveLions}, \cite{Weinstein1,Weinstein2} or
\cite{ShatahStrauss,GrillakisShatahStrauss}) give that any solution
stays near the PLS manifold. Unfortunately, orbital stability
criteria applied to Schrödinger equations use the fact that
Schrödinger operators are bounded from below. Hence the question of
orbital stability for Dirac standing waves cannot be solved by a
straightforward application of the methods used in the Schrödinger
case.

Concerning the asymptotic stability, in the Schrödinger equation,
the question has been solved in several cases. For small stationary
solutions in the simple eigenvalue case it has been studied by
Soffer and Weinstein \cite{SofferWeinstein,SofferWeinstein2}, Pillet
and Wayne \cite{PilletWayne} or Gustafson, Nakanishi and Tsai
\cite{GustafsonNakanishiTsai}. For the two eigenvalue case under a
resonance condition for an excited state, the problem has been
studied by Tsai and Yau
\cite{TsaiYau,TsaiYau2,TsaiYau3,TsaiYau4,Tsai} or Soffer and
Weinstein \cite{SofferWeinstein4,SofferWeinstein5}. Another problem
has been studied by Cuccagna
\cite{Cuccagna2,Cuccagna3,Cuccagna3Err}, he considered the case of
big PLS, when the linearized operator has only one eigenvalue and
obtained the asymptotic stability of the manifold of ground states.
Schlag \cite{Schlag} proved that any ground state of the cubic
nonlinear Schrödinger equation in dimension $3$ is orbitally
unstable but posseses a stable manifold of codimension $9$.


We also would like to mention the works of Buslaev and Perel'mann
\cite{BuslaevPerelman,BuslaevPerelman2,BuslaevPerelman3,BuslaevPerelman4},
Buslaev and Sulem~\cite{BuslaevSulem,BuslaevSulem2}, Weder
\cite{Weder} or Krieger and Schlag \cite{KriegerSchlag} in the one
dimensional Schr\"odinger case. Krieger and Schlag
\cite{KriegerSchlag} proved a result similar to \cite{Schlag} in the
one dimensional case.

In \cite{NBStableDirectionSmallSolitonNR}, we prove that there are
stable directions for the PLS manifold under a  non resonance
assumption on the spectrum of $H:=D_m+V$. This gives a stable
manifold, containing the PLS manifold. But we were not able to say
anything about solutions starting outside the stable manifold.

The results we present here state the existence of a stable
manifold and describe the behavior of solutions starting outside of
it. In fact, we prove the instability of the stable manifold. We
also prove stabilization towards stationary solutions inside the
stable manifold for $H^{s}$ perturbation with $s>2$. We have been
able to obtain it since we impose a resonance condition (see
Assumption \ref{assumption:Resonance} below), while in
\cite{NBStableDirectionSmallSolitonNR}, we assumed there is no
resonance phenomena.

When the perturbations are localized, we can push further this study and 
we obtain a nonlinear scattering.

\medskip

This paper is organized as follow.

In Section \ref{Sec:MainResult}, we present our main results and the
assumptions we need. Subsection \ref{Sec:DecayEstimates}, is devoted
to the statement of the time decay estimates of the propagator
associated with $H=D_m+V$ on the continuous subspace. One is a kind
of smoothness result, in the sense of Kato (see e.g. \cite{Kato}),
the other is a Strichartz type result. We prove these estimates with
the propagation and dispersive estimates proved in
\cite{NBStableDirectionSmallSolitonNR}. In subsection
\ref{Sec:PLSManifold}, we state the existence of small stationary
states forming a manifold tangent to an eigenspace of $H$. The study
of the dynamics around such states leads us to our main results, see
Subsection \ref{Sec:Stability} and \ref{Sec:Scattering}. In 
Subsection \ref{Sec:Stability}, we split a neighborhood of a
stationary state in different parts, each one giving rise to
stabilization or instability. In Subsection \ref{Sec:Scattering}, we state 
our scattering result.

To prove our theorems, we consider our nonlinear system as a
small perturbation of a linear equation. More precisely in
Subsection \ref{Sec:CSUMan}, we show that the spectral properties of
the linearized operator around a stationary state, presented in
Section \ref{Sec:LinearizedOp}, permits to obtain, like in the
linear case, some properties of the dynamics around a stationary
state. We obtain center, center-stable and center-unstable
manifolds. In Section \ref{Sec:Stabilization}, we obtain, with our
time decay estimates, a stabilization towards the PLS manifold for
$H^{s}$ perturbation with $s>2$ in the center manifold. Section
\ref{Sec:Outside} deals with the dynamic outside the center
manifold. Eventually in Section \ref{Sec:EndProof}, we conclude our
study.

Our results are the analogue, in the
Dirac case, of some results of Tsai and Yau \cite{TsaiYau3}, Soffer
and Weinstein \cite{SofferWeinstein}, Pillet and Wayne
\cite{PilletWayne} and Gustafson, Nakanishi and Tsai
\cite{GustafsonNakanishiTsai} about the semilinear Schr\"odinger
equation.


\section{Assumptions and statements}
                                                                    \label{Sec:MainResult}

\subsection{Time decay estimates}
                                                                    \label{Sec:DecayEstimates}
We generalize to small nonlinear perturbations, stability results
for linear systems. These results, like in
\cite{NBStableDirectionSmallSolitonNR}, follow from linear decay
estimates. Here we use smoothness type and Strichartz type estimates
deduced from propagation and dispersive estimates of
\cite{NBStableDirectionSmallSolitonNR}. Hence, we work within the
same assumptions for $V$ and $D_m+V$:
\begin{assumption}
                                                                \label{assumption:1}
The potential~$V:\R^3\mapsto S_4(\C)$ {\it (self-adjoint~$4\times 4$
matrices)} is a smooth function such that there
exists $\rho>5$ with
\begin{equation*}
\forall \alpha\in\N^3,\,\exists C>0,\,\forall x\in\R^3,\, |\p^\alpha
V|(x)\leq\frac{C}{\langle x\rangle^{\rho+|\alpha|}}.
\end{equation*}
\end{assumption}

We notice that by the Kato-Rellich theorem, the operator
\begin{equation*}
H:=D_m+V
\end{equation*}
is essentially self-adjoint on~${\mathcal C}^\infty_0(\R^3,\C^4)$
and self-adjoint on $H^1(\R^3,\C^4)$.

We also mention that Weyl's theorem gives us that the essential
spectrum of $H$ is $(-\infty,-m]\cup[+m,+\infty)$ and the work of
Berthier and Georgescu \cite[Theorem 6, Theorem
A]{BerthierGeorgescu}, gives us that there is no embedded
eigenvalue. Hence the thresholds $\pm m$ are the only points of the
continuous spectrum which can be associated with wave of zero
velocity. These waves perturb the spectral density and diminish the
decay rate in the propagation and the dispersive estimates. We will 
work (like in \cite{NBStableDirectionSmallSolitonNR}) within the 
\begin{assumption}
                                                                \label{assumption:2}
The operator~$H$ presents no resonance at thresholds and no
eigenvalue at thresholds.
\end{assumption}
A resonance is a stationary solution in $H^{1/2}_{-\sigma}\setminus
H^{1/2}$ for any $\sigma \in (1/2,\rho-2)$, where $H^{^t}_{\sigma}$
is given by
\begin{definition}[Weighted Sobolev space]
The weighted Sobolev space is defined by
\begin{equation*}
H_{\sigma}^{t}(\R^3,\C^4)=\left\{ f\in {\mathcal S}'(\R^3),\,
\|\langle Q\rangle^\sigma\langle P\rangle^t f\|_2 <\infty\right\}
\end{equation*}
for~$\sigma,t\in\R$. We endow it with the norm
\begin{equation*}
\|f\|_{ H_{\sigma}^{t}}= \|\langle Q\rangle^\sigma\langle P\rangle^t
f\|_2.
\end{equation*}
If~$t=0$, we write~$L^2_\sigma$ instead of~$H_{\sigma}^{0}$.
\end{definition}
We have used the usual notations~$\langle u\rangle=\sqrt{1+u^2}$,
$P=-i\nabla$, and~$Q$ is the operator of multiplication by~$x$ in
$\R^3$.

\medskip

Now let
\begin{equation*}
{\mathbf P}_c(H)=\one_{(-\infty,-m|\cup[+m,+\infty)}(H)
\end{equation*}
be the projector associated with the continuous spectrum of~$H$ and
${\mathcal H}_c$ be its range. Using \cite[Theorem
1.1]{NBStableDirectionSmallSolitonNR}, we obtain a Limiting
Absorption Principle which gives the $H$-smoothness of 
$\left\langle Q\right\rangle^{-1}$ in the sense of Kato:
\begin{theorem}[Kato smoothness estimates]
                                                                \label{Thm:Smoothness}
If Assumptions~\ref{assumption:1} and~\ref{assumption:2} hold. Then
for any $\sigma\geq 1$ and $s\in \R$, one has:
\begin{align*}
\left\|\left\langle Q\right\rangle^{-\sigma}e^{-\rmi t H}{\mathbf
P}_c\left(H\right)\psi\right\|_{L_t^ 2(\R,H^s(\R^3,\C^4))}&\leq
C\left\|\psi\right\|_{H^s(\R^3,\C^4)},\tag{i}\label{Estimate:Smoothness1}\\
\left\|\int_{\R}e^{\rmi t H}{\mathbf P}_c\left(H\right)\left\langle
Q\right\rangle^{-\sigma}F(t)\;dt\right\|_{H^s(\R^3,\C^4)}&\leq
C\left\|F\right\|_{L_t^2(\R,H^s)},\tag{ii}\label{Estimate:Smoothness2}\\
\left\|\int_{s<t}\left\langle Q\right\rangle^{-\sigma}e^{-\rmi(t-s)
H}{\mathbf P}_c\left(H\right)\left\langle
Q\right\rangle^{-\sigma}F(s)\;ds\right\|_{L_t^2(\R,H^s(\R^3,\C^4))}&\leq
C\left\|F\right\|_{L_t^2(\R,H^s(\R^3,\C^4))}.\tag{iii}\label{Estimate:Smoothness3}\\
\end{align*}
\end{theorem}
\begin{proof}
We first prove \eqref{Estimate:Smoothness1}. For $s=0$, it is (see
{\it e.g.} \cite[Proposition 7.11]{AmreinBoutetdeMonvelGeorgescu} or
\cite[Theorem XIII.25]{ReedSimon4}) a consequence of the limiting
absorption principle:
\begin{equation}\label{Identity:LAP}
\sup_{\Im z\in(0,1)}\left\{\left\|\left\langle
Q\right\rangle^{-\sigma}\left(H-z\right)^{-1}P_c(H)\left\langle
Q\right\rangle^{-\sigma}\right\|_2\right\}<\infty
\end{equation}
which follows from \cite[Theorem
1.1]{NBStableDirectionSmallSolitonNR} or (Theorem \ref{Thm:Propagation} below) 
for $\sigma>5/2$ using the 
fact that the Fourier transform in time of the propagator is the 
resolvent. Actually, the Fourier transform of 
$$\left\langle Q\right\rangle^{-\sigma} e^{-\rmi t(H-\rmi\e)}{\mathbf
P}_c\left(H\right)\bone_{\R^*_+}(t)\left\langle
Q\right\rangle^{-\sigma}f$$ 
in time is 
$$\left\langle
Q\right\rangle^{-\sigma} (H-\lambda-\rmi\e)^{-1}{\mathbf
P}_c\left(H\right)\left\langle Q\right\rangle^{-\sigma}f$$ 
for $f \in L^2(\R^3,\C^4)$. Then we use Born expansion
\begin{equation*}
(H-z)^{-1}=(D_m-z)^{-1}-(D_m-z)^{-1}V(D_m-z)^{-1}
+(D_m-z)^{-1}V(H-z)^{-1}V(D_m-z)^{-1}
\end{equation*}
the limiting absorption in ~\cite[Theorem 2.1(i)]{IftimoviciMantoiu}
(they prove the identity \eqref{Identity:LAP} for $H=D_m$ when
$\sigma=1$) and the fact that
\begin{equation*}
\left\|(H-z)^{-1}\left(1-P_c(H)\right)\right\|_{{\mathcal
B}(L^2)}\leq
\frac{1}{\ds\inf_{\lambda\in(-\infty,-m]\cup[+m,+\infty)}\left|z-\lambda\right|}
\end{equation*}
to obtain \eqref{Identity:LAP} for $\sigma=1$. Hence we have
concluded the proof for $s=0$ and $\sigma\geq1$. For $s\in 2\Z$ and
$\sigma\geq 1$ it follows from the previous cases using boundedness
of $<H>^s<D_m>^{-s}$ and $<H>^{-s}<D_m>^{s}$ (which follow from the
boundedness of $V$ and its derivatives) and the boundedness of
$\langle Q\rangle^ {\mp\sigma}[\langle
Q\rangle^{\pm\sigma},<H>^{s}]\langle H\rangle^ {-s}$ (which follow
from multicommutator estimates see \cite[Appendix
B]{HunzikerSigal}). The rest of the claim \eqref{Estimate:Smoothness1} 
follows by interpolation.

Estimates \eqref{Estimate:Smoothness1} and
\eqref{Estimate:Smoothness2} are equivalent by duality.

To prove estimate \eqref{Estimate:Smoothness3} when $s=0$ (the
general case will follow by the same way as above), we notice 
%
that we have to prove that there exists $C>0$ such that for all
$F,G\in L^2_t(\R,L^2(\R^3,\C^4))$, we have
\begin{multline*}
\left|\iint_{\R^2}\left\langle G(t),\left\langle
Q\right\rangle^{-\sigma} e^{-\rmi (t-s)H}{\mathbf
P}_c\left(H\right)\bone_{\R^*_+}(t-s)\left\langle
Q\right\rangle^{-\sigma}F(s)\right\rangle\;dsdt \right|\\\leq
C\left\|G\right\|_{L_\lambda^2(\R,L^2(\R^3,\C^4))}\left\|F\right\|_{L_\lambda^2(\R,L^2(\R^3,\C^4))}
\end{multline*}
We can suppose that $F$ and $G$ are smooth functions with compact
support from $\R\times \R^3$ to $\C^4$ and we just need to
prove that there exists $C>0$  such that for all $\e>0$,
 for all $F,G\in {\mathcal C}^\infty_0(\R,L^2(\R^3,\C^4))$,
we have
\begin{multline*}
\left|\iint_{\R^2}\left\langle G(t),\left\langle
Q\right\rangle^{-\sigma} e^{-\rmi (t-s)(H-\rmi \e)}{\mathbf
P}_c\left(H\right)\bone_{\R^*_+}(t-s)\left\langle
Q\right\rangle^{-\sigma}F(s)\right\rangle\;dsdt \right|\\\leq
C\left\|G\right\|_{L_\lambda^2(\R,L^2(\R^3,\C^4))}\left\|F\right\|_{L_\lambda^2(\R,L^2(\R^3,\C^4))}.
\end{multline*}
Then we take the limit as $\e\to 0$ and we will conclude using
density arguments. Let us write $A_\e(t)$ for $\left\langle
Q\right\rangle^{-\sigma} e^{-\rmi t (H-\rmi\e)}{\mathbf
P}_c\left(H\right)\bone_{\R^*_+}(t)\left\langle
Q\right\rangle^{-\sigma}$, we have to prove
\begin{equation*}
\left|\int_{\R}\left\langle G(t),(A_\e* F)(t)\right\rangle\;dt
\right|\leq
C\left\|G\right\|_{L_\lambda^2(\R,L^2(\R^3,\C^4))}\left\|F\right\|_{L_\lambda^2(\R,L^2(\R^3,\C^4))}.
\end{equation*}
Using Plancherel's identity in $L_t^2(\R,L^2(\R^3,\C^4))$ and 
$A_\e*F\in L_t^2(\R,L^2(\R^3,\C^4))$, we
just need to prove
\begin{equation*}
\left|\int_{\R}\left\langle \widehat{G}(\lambda),\widehat{A_\e*
F}(\lambda)\right\rangle\;d\lambda \right|\leq
C\left\|\widehat{G}\right\|_{L_\lambda^2(\R,L^2(\R^3,\C^4))}\left\|\widehat{F}\right\|_{L_\lambda^2(\R,L^2(\R^3,\C^4))}.
\end{equation*}
Since the Fourier transform in time of the propagator is the resolvent, 
$F$ is smooth with compact support and $\e>0$, we obtain  
\begin{equation*}
    \widehat{A_\e*F}(\lambda)=\left\langle
Q\right\rangle^{-\sigma} (H-\lambda-\rmi\e)^{-1}{\mathbf
P}_c\left(H\right)\left\langle
Q\right\rangle^{-\sigma}\widehat{F}(\lambda)
\end{equation*}
Hence we just have to prove
\begin{equation*}
\left|\int_{\R}\left\langle \widehat{G}(\lambda),\left\langle
Q\right\rangle^{-\sigma} (H-\lambda-\rmi\e)^{-1}{\mathbf
P}_c\left(H\right)\left\langle
Q\right\rangle^{-\sigma}\widehat{F}(\lambda) \right\rangle\;d\lambda
\right|\leq
C\left\|\widehat{G}\right\|_{L_\lambda^2(\R,L^2(\R^3,\C^4))}\left\|\widehat{F}\right\|_{L_\lambda^2(\R,L^2(\R^3,\C^4))}.
\end{equation*}
This in turn follows from the Limiting Absorption Principle
\eqref{Identity:LAP} just proved.
\end{proof}

To state the next result, we need the
\begin{definition}[Besov space]
For~$s\in\R$ and~$1\leq p,q\leq\infty$, the Besov
space~$B^s_{p,q}(\R^3,\C^4)$ is the space of all~$f\in {\mathcal
S}'(\R^3,\C^4)$ (dual of the Schwartz space) such that
\begin{equation*}
\ds \|f\|_{B^s_{p,q}}= \left(\sum_{j\in\N}2^{jsq} \|\varphi_j *
f\|_p^q\right)^{\frac{1}{q}}<+\infty
\end{equation*}
with~$\widehat{\varphi} \in {\mathcal C}^\infty_0(\R^n\setminus
\left\{0\right\})$ such that~$\sum_{j\in
\Z}\widehat{\varphi}(2^{-j}\xi)=1$ for all~$\xi
\in\R^3\setminus\left\{0\right\}$,
$\widehat{\varphi}_j(\xi)=\widehat{\varphi}(2^{-j} \xi)$ for all
$j\in\N^*$ and for all~$\xi \in\R^3$, and
$\widehat{\varphi_0}=1-\sum_{j\in \N^*}\widehat{\varphi}_j$. It is
endowed with the natural norm~$f\in B^s_{p,q}(\R^3,\C^4) \mapsto
\|f\|_{B^s_{p,q}}$.
\end{definition}
Using the Dispersive estimates of \cite[Theorem
1.2]{NBStableDirectionSmallSolitonNR} and \cite[Theorem
10.1]{KeelTao}, we obtain the
\begin{theorem}[Strichartz estimates]
                                                                \label{Thm:Strichartz}
If Assumptions \ref{assumption:1} and \ref{assumption:2} hold. Then
for any $2\leq p,q \leq \infty$, $\theta\in[0,1]$, with
$(1-\frac{2}{q})(1\pm\frac{\theta}{2})=\frac{2}{p}$ and
$(p,\theta)\neq(2,0)$, and for any reals $s$, $s'$ with $s'-s\geq
\alpha(q)$ where $\alpha(q)=(1+\frac{\theta}{2})(1-\frac{2}{q})$,
there exists a positive constant $C$ such that
\begin{align}
\|e^{-\rmi t  H}P_c(H)\psi\|_{L_t^{p}(\R,B^s_{q,2}(\R^3,\C^4))}&\leq
    C\|\psi\|_{H^{s'}(\R^3,\C^4)},
    \tag{i}\label{Eq:Strichartz_1}\\
\|\int e^{\rmi t H}P_c(H) F(t)\,dt\|_{H^s}&\leq
    C\|F\|_{L_t^{p'}(\R,B^{s'}_{q',2}(\R^3,\C^4))},\tag{ii}\label{Eq:Strichartz_2}\\
\| \int_{s<t} e^{-\rmi(t-s)H}P_c(H)
F(s)\,ds\|_{L_t^{p}(\R,B^{-s}_{q,2}(\R^3,\C^4))}&\leq C
    \|F\|_{L_t^{\tilde{p}'}(\R,B^{\tilde{s}}_{\tilde{q}',2}(\R^3,\C^4))}\tag{iii}\label{Eq:Strichartz_3},
\end{align}
for any $r\in[1,\infty]$, $(\tilde{q}, \tilde{p})$ chosen like $(q,
p)$ and $s+\tilde{s}\geq\alpha(q)+\alpha(\tilde{q})$.
\end{theorem}
\begin{proof}
This is a consequence of \cite[Theorem 10.1]{KeelTao} applied to
$U(t)=e^{-\rmi t H}P_c(H)$, using \cite[Theorem
1.2]{NBStableDirectionSmallSolitonNR} or Theorem \ref{Thm:Dispersion} below
and
\begin{equation*}
B^{(1+\frac{\theta}{2})(1-\frac{2}{q})+s}_{q,2}\hookrightarrow(H^s,
B^{1+\theta/2+s}_{1,2})_{2/((1\pm \theta/2)p),2}
\end{equation*}
continuously for $p\geq 2$ ($p\neq 2$ if $\theta=0$) 
and $1/q=1-1/((1\pm \theta/2)p)$. For these embeddings, we refer to the
 proof of \cite[Theorem 6.4.5]{BerghLofstrom} as well as the properties of the real
interpolation (see \cite{BerghLofstrom} or \cite{Triebel}). More
precisely for $\theta=0$ or $1$ it is obvious. In the other cases, we work like
in proof of \cite[Theorem 6.4.5(3)]{BerghLofstrom}:

We use \cite[Theorem 6.4.3]{BerghLofstrom} ($B^s_{p,2}$ is a
retract of $l^s_2(L^p)$ for $s\in \R$ and $p,q\in[1,\infty]$) and
\cite[Theorem 5.6.2]{BerghLofstrom} (about the interpolation of
$l^s_2(L^p)$ spaces) with \cite[Theorem
5.2.1]{BerghLofstrom} (about the interpolation of $L^p$ spaces). Then
we conclude using the injection of $L^p$ spaces into some Lorentz
spaces \cite[Section 1.3 \& Exercice 1.6.8]{BerghLofstrom}.

In the case $\theta\neq 0$, the proof is actually simpler. We can
prove it using the usual $TT^*$ method and the H\"older inequality
instead of the Hardy-Littlewood-Sobolev inequality.

\end{proof}

\subsection{The manifold of PLS}
                                                                    \label{Sec:PLSManifold}
We study the following nonlinear Dirac equation
\begin{equation}
\left\{\begin{array}{l} \rmi\p_t \psi=H\psi+\nabla F(\psi)\\
\psi(0,\cdot)=\psi_0.
\end{array}\right.
                                                                \label{Eq:NLD}
\end{equation}
with~$\psi\in {\mathcal C}^1(I,H^1(\R^3,\C^4))$ for some open
interval~$I$ which contains~$0$ and~$H=D_m+V$. The nonlinearity~$F:
\C^4\mapsto \R$ is a differentiable map for the real structure
of~$\C^4$ and hence the~$\nabla$ symbol has to be understood for the
real structure of~$\C^4$. For the usual hermitian product of~$\C^4$,
one has
\begin{equation*}
DF(v)h=\Re \langle \nabla F(v),h\rangle.
\end{equation*}
If $F$ has a gauge invariance (see Equation
\eqref{Eq:GaugeInvariance} or Assumption
\ref{assumption:NonLinearity} below), this equation may have
stationary solutions {\it i.e.} solution of the form $e^{-\rmi E
t}\phi_0$ where $\phi_0$ satisfies the nonlinear stationary
equation:
\begin{equation*}
    E\phi_0=H\phi_0+\nabla F(\phi_0).
\end{equation*}
We will notice that the Dirac operator $D_m$ have an interesting
invariance property due to its matrix structure. This invariance can
be shared by some perturbed Dirac operators and gives a consequence
of a theorem of Kramers, see \cite{BalslevHelffer,Parisse}. Indeed if
we introduce $K$ the antilinear operator defined by:
\begin{equation}
                                                                \label{Def:OperatorK}
K\left(\begin{array}{c}\psi\\\chi\end{array}\right) =
\left(\begin{array}{c}\sigma_2\overline{\psi}\\
\sigma_2\overline{\chi}\end{array}\right)\mbox{ with }
\sigma_2=\left(\begin{array}{cc}0&-\rmi\\\rmi&0\end{array}\right).
\end{equation}
The operator $D_m$ commutes with $K$. So if $V$ also commutes with
$K$, we obtain that the eigenspaces of $H$ are always of even
dimension. Here we work with the
\begin{assumption}
                                                                \label{assumption:Spectrum}
The potential $V$ commutes to $K$. The operator~$H:=D_m+V$ has only
two double eigenvalues $\lambda_0<\lambda_1$,
with~$\left\{\phi_0,K\phi_0\right\}$
and~$\left\{\phi_1,K\phi_1\right\}$ as associated orthonormalized
basis.
\end{assumption}

We also need the
\begin{assumption}
                                                                \label{assumption:NonLinearity}
The function~$F:\C^4\mapsto \R$ is in~${\mathcal C}^\infty(\R^8,\R)$
and satisfies $F(z)=O(|z|^{4})$ as~$z\to 0$. Moreover, it has the
following invariance properties:
\begin{equation*}
\forall z\in\C^4,\;\forall \theta \in
\R,\;F(Kz)=F(z),\;F(e^{i\theta}z)=F(z).
\end{equation*}
\end{assumption}

We obtain the
\begin{proposition}[PLS manifold]
                                                                \label{Prop:ManifoldPLS}
If Assumptions~\ref{assumption:1}--\ref{assumption:NonLinearity}
hold. Then for any~$\sigma \in \R^+$, there exist~$\Omega$ a
neighborhood of~$0$ in $\C^2$, a smooth map
\begin{equation*}
h: \Omega \mapsto
\left\{\phi_0,K\phi_0\right\}^{{{\bot}}}\cap{H^2(\R^3,\C^4)}\cap
L^2_\sigma(\R^3,\C^4)
\end{equation*}
and a smooth map~$E: \Omega \mapsto \R$ such that
$S((u_1,u_2))=u_1\phi_0+u_2K\phi_0+h((u_1,u_2))$ satisfy for
all~$U\in\Omega$,
\begin{equation}
                                                                \label{Eq:StationaryStates}
H S(U)+\nabla F(S(U))=E(U)S(U),
\end{equation}
with the following properties
\begin{equation*}
\left\{\begin{array}{l}
h((u_1,u_2))=\left(\frac{u_1}{|(u_1,u_2)|}Id_{\C^4}+\frac{u_2}{|(u_1,u_2)|}
K\right)h\left((|(u_1,u_2)|,0)\right),\quad\forall U=(u_1,u_2) \in \Omega,\\
h(U)=O(|U|^2),\\
E(U)=E(\left|U\right|),\\
E(U)=\lambda_0+O(|U|^2).
\end{array}\right.
\end{equation*}
\end{proposition}
\begin{proof}
This result is adapted from \cite[Proposition 2.2]{PilletWayne}
after the reduction due to the invariance of the problem with
respect to $K$.
\end{proof}

Moreover, we have
\begin{lemma}[exponential decay]
                                                                \label{lemma:ExpDecay}
For any~$\beta\in \N^4$,~$s\in \R^+$ and~$p,q\in[1,\infty]$. There
exist~$\gamma
>0$,~$\e>0$ and~$C>0$ such that for all~$U \in B_{\C^2}(0,\e)$ one has
\begin{equation*}
\|e^{\gamma \langle Q\rangle}\p_U^\beta S(U)\|_{B^s_{p,q}} \leq C
\|S(U)\|_2,
\end{equation*}
where $\p_{(u_1,u_2)}^\beta=\frac{\p^{|\beta|}}{\p^{\beta_1}\Re u_1
\p^{\beta_2}\Im u_1\p^{\beta_3}\Re u_2 \p^{\beta_4}\Im u_2}$.
\end{lemma}
\begin{proof}
This is proved like in \cite[Lemma
4.1]{NBStableDirectionSmallSolitonNR}, where we used ideas of
\cite{Hislop}.
\end{proof}

\subsection{The unstable manifold and the stabilization}
                                                                    \label{Sec:Stability}
Each stationary solution previously introduced has, like in
\cite{NBStableDirectionSmallSolitonNR}, a stable manifold. Under the
following assumption, we can prove that the stable manifold is
unstable, that is to say that a small perturbation of a stationary
solution starting outside of this manifold leaves any neighborhood
of this stationary solution. We work with the
\begin{assumption}
                                                                \label{assumption:Resonance}
The  resonant condition
\begin{equation*}
|\lambda_1-\lambda_0|>\min\{|\lambda_0+m|,\,|\lambda_0-m|\}
\end{equation*}
holds. Moreover, we have the Fermi Golden Rule
\begin{equation}\label{FermiGoldenRule}
    \Gamma(\phi)=\lim_{\substack{\e\to 0,\\\e>0}}\left\langle d^2
    F(\phi)\phi_1,\Im\left((H-\lambda_0)
    +(\lambda_1-\lambda_0)-\rmi\e\right)^{-1}P_c(H)d^2
    F(\phi)\phi_1\right\rangle>0
\end{equation}
for any non zero eigenvector $\phi$  associated with $\lambda_0$.
\end{assumption}
In this assumption, the notation $d^2F$ denotes the differential of
$\nabla F$ with respect to the real structure of $\C^4$.

Let us introduce the linearized operator $JH(U)$ around a stationary
state $S(U)$:
\begin{equation*}
    H(U)=H+d^2F(S(U))-E(U).
\end{equation*}
We notice that the operator $H(U)$ is not
$\C$-linear but only $\R$-linear. Hence we work with the space $L^2(\R^3,\R^4\times
\R^4)$ instead of $L^2(\R^3,\C^4)$ by writing
\begin{equation*}
    \begin{pmatrix}\Re \phi\\ \Im \phi\end{pmatrix}
\end{equation*}
instead of $\phi$. The multiplication by $-\rmi$ becomes the operator
\begin{equation*}
    J=\begin{pmatrix}0 &-I_{\R^4}\\
    I_{\R^4}&0\\\end{pmatrix}.
\end{equation*}

Now we mention some spectral properties of the \emph{real operator}
$JH(U)$ in $L^2(\R^3,\C^4\times \C^4)$ (the
complexified of $L^2(\R^3,\R^4\times \R^4)$) which are needed to state and
to understand our main theorem. These properties will be proven in
subsection \ref{Sec:LinearizedOp}.

\begin{proposition}[Spectrum of $JH(U)$]
The operator $JH(U)$ in $L^2(\R^3,\C^4\times \C^4)$  has a four dimensional
geometric kernel and four double eigenvalues $E_1(U)$,
$\overline{E_1}(U)$, $-E_1(U)$ and $-\overline{E_1}(U)$ with $\Re
E_1(U)>0$.

The eigenspaces associated with $E_1(U)$ and $\overline{E_1}(U)$ are
conjugated via the complex conjugation.
The same holds for $-E_1(U)$ and $-\overline{E_1}(U)$.

The rest of the spectrum is the essential (or continuous) spectrum.
We write ${\mathcal H}_c(U)$ for the space associated with the
continuous spectrum. The space ${\mathcal H}_c(U)$ is the orthogonal
of the previous eigenspaces and the geometric kernel of $JH(U)$ and is 
invariant by the complex conjugation.
\end{proposition}
\begin{proof}
See subsection \ref{Sec:LinearizedOp} below.
\end{proof}
We will work on the real part of the sum the eigenspaces 
associated with $E_1(U)$ and $\overline{E_1}(U)$: 
$X_u(U)\subset L^2(\R^3,\R^4\times \R^4)$, we
introduce a real basis $\left(\xi_i(U)\right)_{i=1,\ldots,4}$ of
$X_u(U)$.
We will also work in
the real part of the sum of the eigenspaces associated with $-E_1(U)$ 
and $-\overline{E_1}(U)$ :
 $X_s(U)\subset L^2(\R^3,\R^4\times \R^4)$, we introduce a real basis
$\left(\xi_i(U)\right)_{i=5,\ldots,8}$ of $X_s(U)$.

We can state our main theorems which will be proved in the 
sections \ref{Sec:LinearizedOp}, \ref{Sec:Stabilization}, 
\ref{Sec:Outside} and \ref{Sec:EndProof}.
\begin{theorem}[Central manifold and asymptotic stability]      \label{Thm:Main}
If Assumptions \ref{assumption:1}--\ref{assumption:Resonance} hold.
Then, for $s>\beta+2>2$ and
 $\sigma>3/2$, there exist $\e>0$, a continuous map
 $r : B_{\C^2}(0,\e)\mapsto \R$ with $r(U)=O(\left|U\right|^2)$, $C>0$,
${\mathcal V}$ a neighborhood of $(0,0)$ in \begin{equation*}
    {\mathcal S}=\left\{\left(U,z\right);\;U\in B_{\C^2}(0,\e),\; z\in {\mathcal H}_c(U)\cap B_{H^{s}}(0,r(U)) \right\}
\end{equation*}
endowed with the metric of $\C^2\times H^{s}$ and a  map $\Psi: {\mathcal V}
  \mapsto \R^8$, smooth on ${\mathcal V}\setminus{(0,0)}$
  satisfying  for any non zero $U\in B_{\C^2}(0,\e)$
  $$\left\|\Psi(U,z)\right\|=O(\left\|z\right\|_{H^s}^2)$$
for all $z\in{\mathcal H}_c(U)\cap B_{H^{s}}(0,r(U))$ with $(U,z)\in
{\mathcal V}$ such that the following is true.

For any initial condition of the form
\begin{equation*}
    \psi_0=S(U_0)+z_0+A\cdot\xi(U_0)
\end{equation*}
with $(U_0,z_0)\in {\mathcal V}$ and  $A=\Psi(U_0,z_0)$, there exists
a solution $\psi\in \cap_{k=0}^2 {\mathcal
C}^k\left(\R,H^{s-k}\right) $ of \eqref{Eq:NLD} with initial
condition $\psi_0$ and this solution is unique in
$L^\infty((-T,T),H^{s}(\R^3,\C^4))$ for any $T>0$.

Moreover, we have for all $t\in\R$
\begin{equation}                                                \label{eq:Stabilization}
    \psi(t)=e^{-\rmi \int_0^t E(U(v))\;dv}S(U(t))+\e(t)
\end{equation}
with $\left\|\dot{U}\right\|_{L^q(\R)}\leq C\|z_0\|_{H^{s}}^2$ for
all $q\in[1,\infty]$, $\ds\lim_{t\to\pm \infty} U(t)=U_{\pm\infty}$
and
\begin{equation*}
\max\Big\{\left\|\e\right\|_{L^\infty(\R^\pm,
H^{s})},\;\left\|\e \right\|_{L^2(\R^\pm, H^{s}_{-\sigma})},
\left\|\e \right\|_{L^2(\R^\pm, B^\beta_{\infty,2})}\Big\}\leq
\ds C\|z_0\|_{H^{s}}.
\end{equation*}
\end{theorem}

\begin{theorem}[Center stable end center unstable manifold]     \label{Thm:MainOutside}
With the same assumptions and notations as Theorem \ref{Thm:Main}, 
let $\mathcal{CM}$ be the graph of $(U,z)\in {\mathcal V}\mapsto 
S(U)+z +\Psi(U,z)\cdot\xi(U)$ then for the set
\begin{equation*}
    \widetilde{\mathcal S}=\left\{\left(U,z,p\right);\;U\in B_{\C^2}(0,\e),\; z\in {\mathcal H}_c(U)\cap B_{H^{s}}(0,r(U)),
    p\in B_{\R^4}(0,r(U)),\right\}
\end{equation*}
endowed with the metric of $\C^2\times H^{s}\times\R^4$, there exist $C>0$, neighborhoods ${\mathcal W}_\pm$ of $(0,0,0)$ in $\widetilde{\mathcal S}$ and
maps
  $\Phi_\pm: {\mathcal W}_\pm\mapsto \R^8$, smooth on ${\mathcal
W}_\pm\setminus\{(0,0,0)\}$ with
$$\left\|\Phi_\pm(U,z,p)\right\|=O(\left\|z\right\|_{H^s}^2+\|p\|^2)$$
for all $(U,z,p)\in {\mathcal W}_\pm$ such that for any initial condition of the form
\begin{equation*}
\psi_0=S(U_0)+z_0+(p_+,p_-).\xi(U_0)
\end{equation*}
not in $\mathcal{CM}$, the following is
true.
\begin{enumerate}
	\item If $(U_0,z_0, p_+)\in {\mathcal W}_+$ and 
	$p_-= \Phi_+(U_0,z_0, p_+)$ ({\it resp.} If $(U_0,z_0, p_-)\in {\mathcal W}_-$ and 
	$p_+= \Phi_-(U_0,z_0, p_-)$) then for any small neighborhood ${\mathcal O}$ of $S(U_0)$
containing $\psi_0$ there exist $t_\pm(\psi_0)>0$ and  a solution
$\psi_+\in\cap_{k=0}^2{\mathcal C}^k([-t_+; +\infty),H^{s-k})$ (resp.
$\psi_-\in\cap_{k=0}^2{\mathcal C}^k((-\infty;t_-],H^{s-k})$)  of
\eqref{Eq:NLD} with initial condition $\psi_0$ and this solution is
unique in $L^\infty((-T',T),H^{s}(\R^3,\C^4))$ for any $T>0$ (resp
$T\in(0,t_-)$) and any $T'\in (0,t_+)$ (resp $T'<0$).

Moreover, there exist $C>0$, $\phi_\pm(t)\in \mathcal{CM}$ and
$\rho_+(t)\in X_s(U_0)$ (resp. $\rho_-(t)\in X_u(U_0))$) for all
$t>-t_-$ (resp for all $t<t_+$) such that $\psi_\pm(t)=\phi_\pm(t)+\rho_\pm(t)$
with
\begin{equation*}
    \left\|\rho_\pm(t)\right\|_{H^{s}}\leq C \left\|\rho_\pm(0)\right\|_{H^{s}}e^{\mp\gamma t}
    \mbox{ as }t\to \pm\infty\quad\mbox{ and }\quad
     \psi_\pm(\mp t_\pm)\notin {\mathcal O}
\end{equation*}
where $\gamma$ is in a ball around $1/2 \Gamma(U_0)$, the radius of
which is $O(|U_0|^6))$.

We also have
\begin{equation*}
    \phi_\pm(t)=e^{-\rmi \int_0^t E(U_\pm(v))\;dv}S(U_\pm(t))+\e_\pm(t), \quad
\forall t>t- \left (\mbox{resp. }  \forall t<t+\right)
\end{equation*}
with $\left\|\!\dot{U}_+ \!\right\|_{L^q((-t_-,+\infty))}\!\!\leq\!
C\!\left(\|z_0\|_{H^{s}}\!\!+\!\|\rho_\pm(0)\|_{H^{s}}\right)^2$ (resp.
$\left\|\!\dot{U}_- \!\right\|_{L^q((-\infty,t_+))}\!\!\leq\!
C\!\left(\|z_0\|_{H^{s}}\!\!+\!\|\rho_\pm(0)\|_{H^{s}}\right)^2$) for all $q\in[1,\infty]$, $\ds\lim_{t\to\pm
\infty} U_\pm(t)=U_{\pm\infty}$ and
\begin{multline*}
\max\Big\{\left\|\e_+\right\|_{L^\infty((-t_-,+\infty),
H^{s})},\;\left\|\e_+ \right\|_{L^2((-t_-,+\infty),
H^{s}_{-\sigma})}, \left\|\e_+ \right\|_{L^2((-t_-,+\infty),
B^\beta_{\infty,2})}\Big\}\\\leq \ds
C\left(\|z_0\|_{H^{s}}+\|\rho_\pm(0)\|_{H^{s}}\right),
\end{multline*}
\begin{multline*} \Bigg(\mbox{resp. }
\max\Big\{\left\|\e_-\right\|_{L^\infty((-\infty,t_+),
H^{s})},\;\left\|\e_{-} \right\|_{L^2((-\infty,t_+),
H^{s}_{-\sigma})}, \left\|\e_{-} \right\|_{L^2((-\infty,t_+),
B^\beta_{\infty,2})}\Big\}\\\leq \ds
C\left(\|z_0\|_{H^{s}}+\|\rho_\pm(0)\|_{H^{s}}\right)\Bigg).
\end{multline*}

\item If $(U_0,z_0, p_+)\in {\mathcal W}_+$ and 
	$p_-\neq  \Phi_+(U_0,z_0, p_+)$ or $(U_0,z_0, p_-)\in {\mathcal W}_-$ and 
	$p_+\neq \Phi_-(U_0,z_0, p_-)$, then there exist $t_+(\psi_0)>0$, $t_-(\psi_0)<0$
and a unique solution $\psi$ of \eqref{Eq:NLD} with initial
condition $\psi_0$ such that
    for any small neighborhood ${\mathcal O}$ of $S(U_0)$
    containing $\psi_0$, $\phi\in \cap_{k=0}^2{\mathcal C}^k([t_-;t_+],H^{s-k})$ with
    $\psi(t_+)\notin {\mathcal
    O}$ and $\psi(t_-)\notin {\mathcal O}$. This solution is unique in
$L^\infty((T',T),H^{s}(\R^3,\C^4))$ for any $T\in(0,t_+)$ and any
$T'\in (t_-,0)$.
\end{enumerate}
\end{theorem}

\medskip

The first theorem shows, as in
\cite{NBStableDirectionSmallSolitonNR}, that perturbations in the
direction of the continuous subspace, except four directions, relax
towards stationary solutions. 
We have excluded four directions in the continuous subspace, which,
due to resonance phenomena, induce orbital instability. The second
theorem tells us what happens for perturbations in the
directions of an excited state and in the four directions of the
continuous spectrum for which we haven't the stabilization. We
thus study eight directions: four of them give a manifold on which
there hold exponential stabilization in positive time and orbital
instability in negative time, while the four others give a manifold
on which there hold exponential stabilization in negative time and
orbital instability in positive time. Outside these manifolds, we
have orbital instability in both negative and positive time.

\bigskip

\begin{remark}

In the Theorem, we notice that when $U_0=0$ then $z_0=0$ and $p=0$
so the theorem do not say anything for this case. In fact, the
charge conservation gives the orbital stability of $0$. But we
cannot extend the previous results to $0$ since we can build a manifold of
stationary states tangent to the eigenspace associated with
$\lambda_1$ similarly to Proposition \ref{Prop:ManifoldPLS}.

\end{remark}

\subsection{The nonlinear scattering}
                                                                    \label{Sec:Scattering}
If we choose a localized $z_0$, we are able to expand further 
\eqref{eq:Stabilization} as stated by the following theorems also 
proved in sections \ref{Sec:LinearizedOp}, \ref{Sec:Stabilization}, 
\ref{Sec:Outside} and \ref{Sec:EndProof}. 
\begin{theorem}                                                 \label{Thm:MainScatt}
With the assumptions and the notations of Theorem \ref{Thm:Main}, for the set
\begin{equation*}
    {\mathcal S}_\sigma=\left\{\left(U,z\right);\;U\in B_{\C^2}(0,\e),\; z\in {\mathcal H}_c(U)\cap B_{H^{s}_\sigma}(0,r(U)) \right\}
\end{equation*}
endowed with the metric of $\C^2\times H^{s}_\sigma$, there exists a neigborhhood $\mathcal{V}_\sigma$ of $(0,0)$ in ${\mathcal S}_\sigma$ such that the following is true. If  $A=\Psi(U_0,z_0)$ with $(U_0,z_0)\in
\mathcal{V}_\sigma$,
 there exist $\mathcal{V}^\pm_\sigma$ open neighborhoods of $(0,0)$ in ${\mathcal S}_\sigma$ and $\left(U_{\pm\infty};z_{\pm\infty}\right)\in\mathcal{V}^\pm_\sigma$,
such that
\begin{equation*}
\left|V_{ \pm\infty}-U_0\right|\leq C\|z_0\|_{H_\sigma^{s}}^2,\; \left\|z_{
\pm\infty}-z_0\right\|_{H^{s}}\leq
C\|z_0\|_{H^s_\sigma}^2,
\end{equation*}
and for all $t\in\R$
\begin{equation*}
    \psi(t)=e^{-\rmi t E(V_{\pm\infty})}S(V_\pm(t))+e^{J t E(V_{\pm\infty})}e^{J t H(V_{\pm\infty})}z_{\pm\infty}+\e_\pm(t)\\
\end{equation*}
with 
\begin{eqnarray*}
\left|\dot{V}_\pm(t)+\rmi\left(E(V_\pm(t))-E(V_{\pm\infty})\right)\right|&\leq& \frac{C}{\left\langle t\right\rangle^2}\|z_0\|_{H^{s}_\sigma}^2,\\ 
\left|V_\pm(t)-V_{\pm\infty}\right|&\leq& \frac{C}{\left\langle t\right\rangle}\|z_0\|_{H^{s}_\sigma},\\ 
\max\left\{\left\|\e_\pm(t)\right\|_{
H^{s}},\;\left\|\e_{\pm}(t) \right\|_{H^{s}_{-\sigma}},
\left\|\e_{\pm}(t) \right\|_{B^\beta_{\infty,2}}\right\}&\leq&
\ds \frac{C}{\left\langle t\right\rangle^2}\|z_0\|_{H^s_\sigma}^2\\
\mbox{and }
\left\|e^{-J t H(V_{\pm\infty})}e^{J \int_0^t (E(V_{\pm}(s))-E(V_{\pm\infty})\;ds}\e_{\pm}(t) \right\|_{H^{s}_{\frac{3}{2}}}&\leq&
\ds \frac{C}{\left\langle t\right\rangle^{\frac{1}{2}}}\|z_0\|_{H^s_\sigma}^2
\end{eqnarray*}
for all $t\in \R$.

Moreover, the maps
\begin{equation*}
\left(U_{0};z_{0}\right)\in\mathcal{V}_\sigma
\mapsto\left(V_{\pm\infty};z_{\pm\infty}\right)\in{\mathcal V}^\pm_\sigma
\end{equation*}
are bijective.


\end{theorem}

\begin{remark}

The fact that $z_0$ is localized gives us the convergence of 
\begin{equation*}
\int_0^tE(U(v))\;dv-tE(U_{\pm\infty})
\end{equation*}
as $t\to \pm \infty$ and allows us to obtain an asymptotic 
profile for the dispersive part of the perturbed solution~$\phi$.  

What we call the nonlinear scattering result is essentially the fact 
that the maps
\begin{equation*}
\left(U_{0};z_{0}\right)\in\mathcal{V}_\sigma
\mapsto\left(U_{\pm\infty};z_{\pm\infty}\right)\in{\mathcal V}^\pm_\sigma,
\end{equation*}
are well defined and bijective (actually the surjectivity is called asymptotic completness).
 
Using wave operators for the couple $(JH(U),JD_m)$, we can obtain an expansion of the form
$\psi(t)=e^{-\rmi\int_0^t E(U(v))\;dv}S(U_{\pm\infty})+e^{-\rmi t
D_m}z_\pm+\e_\pm(t)$ but we will only have 
$$\left\|z_{\pm
\infty}-z_0\right\|_{H^{s}}\leq C\|z_0\|_{H^{s}}$$ 
and
\begin{equation*}
\max\Big\{\left\|\e_\pm \right\|_{L^\infty(\R^\pm,
H^{s})},\;\left\|\e_\pm\right\|_{L^2(\R^\pm, H^{s}_{-\sigma})},
\left\|\e_\pm\right\|_{L^2(\R^\pm, B^\beta_{\infty,2})}\Big\}\leq
\ds C\|z_0\|_{H^{s}},
\end{equation*}
or using wave operators for the couple $(JH(U),JH)$ we can 
obtain an expansion  of the form $\psi(t)=e^{-\rmi\int_0^t
E(U(v))\;dv}S(U_{\pm\infty})+e^{-\rmi t H}z_\pm+\e_\pm(t)$ with
$$\left\|z_{\pm \infty}-z_0\right\|_{H^{s}}\leq
C\left(\left|U_0\right|+\|z_0\|_{H^{s}}\right)\|z_0\|_{H^{s}}$$ and
\begin{equation*}
\max\Big\{\left\|\e_\pm \right\|_{L^\infty(\R^\pm,
H^{s})},\;\left\|\e_\pm\right\|_{L^2(\R^\pm, H^{s}_{-\sigma})},
\left\|\e_\pm\right\|_{L^2(\R^\pm, B^\beta_{\infty,2})}\Big\}\leq
\ds C\left(\left|U_0\right|+\|z_0\|_{H^{s}}\right)\|z_0\|_{H^{s}}.
\end{equation*}
But in these cases, we cannot obtain a nice asymptotic 
completness statement.

We have (in the previous expansions)
\begin{multline*}
\max\Big\{\sup_{t\in\R}\left(\left\|e^{J t
H(U_{\pm\infty})}z_{\pm\infty}\right\|_{
H^{s}}\right),\;\sup_{t\in\R}\left(\left\langle t\right\rangle^{3/2}\left\|e^{J t
H(U_{\pm\infty})}z_{\pm\infty}\right\|_{H^{s}_{-\sigma}}\right),\\
 \sup_{t\in\R}\left(\left\langle t\right\rangle^{3/2}\left\|e^{J t
H(U_{\pm\infty})}z_{\pm\infty}\right\|_{B^\beta_{\infty,2}}\right)\Big\}\leq C\|z_{\pm\infty}\|_{H^{s}}.
\end{multline*}
This follows from Lemma \ref{Lem:StabilizationForTloc} and 
Lemma \ref{Lem:OnNonLinearScattering}.
\end{remark}

Outside the center manifold, we can also have an expansion of the
same type. But due to the presence of exponentially stable and
unstable directions, one cannot expect a scattering result of the
same type. Actually we cannot obtain the injectivity of the corresponding mappings. 
We have the
\begin{theorem}                                                 \label{Thm:MainAsymp}
With the assumptions and the notations of
Theorem \ref{Thm:MainOutside}, for the sets
\begin{equation*}
    \widetilde{\mathcal S}_\sigma=\left\{\left(U,z,p\right);\;U\in B_{\C^2}(0,\e),\; z\in {\mathcal H}_c(U)\cap B_{H^{s}_\sigma}(0,r(U)),
    p\in B_{\R^4}(0,r(U)),\right\}
\end{equation*}
endowed with the metric of $\C^2\times H^{s}_\sigma\times\R^4$, there exist $C>0$, neighborhoods ${\mathcal W}^\pm_\sigma$ of $(0,0,0)$ in ${\mathcal S}_\sigma$ such that the following is true.

If $\psi_0\notin
\mathcal{CM}$, $(U_0,z_0, p_+)\in {\mathcal W}^+_\sigma$ and 
	$p_-= \Phi_+(U_0,z_0, p_+)$ ({\it resp.} $(U_0,z_0, p_-)\in {\mathcal W}^-_\sigma$ and 
	$p_+= \Phi_-(U_0,z_0, p_-)$) then there exist $C>0$,  $\phi_\pm(t)\in
\mathcal{CM}$ and $\rho_\pm(t)\in X_s(U_0)$  for all $t>-t_-$ (resp for all $t<t_+$) such that
$\psi_\pm(t)=\phi_\pm(t)+\rho_\pm(t)$ with
\begin{equation*}
    \left\|\rho_\pm(t)\right\|_{H^{s}}\leq C \left\|\rho_\pm(0)\right\|_{H^{s}}e^{\mp\gamma t}
    \mbox{ as }t\to \pm\infty\quad\mbox{ and }\quad
     \psi(\mp t_\pm)\notin {\mathcal O}
\end{equation*}
where $\gamma$ is in a ball around $1/2 \Gamma(U_0)$, the radius of
which is $O(|U_0|^6))$,  there exist
$\left(V_{\pm\infty};z_{\pm\infty}\right)\in {\mathcal S}$ such that
\begin{eqnarray*}
\left|V_{ \pm\infty}-U_0\right|&\leq&
C\left(\|z_0\|_{H^{s}_\sigma}+\|\rho_\pm(0)\|_{H^{s}}\right)^2,\\\; \left\|z_{
\pm\infty}-z_0\right\|_{H^{s}}&\leq&
C\left(\|z_0\|_{H^{s}_\sigma}+\|\rho_\pm(0)\|_{H^{s}}\right)^2,
\end{eqnarray*}
and for all $t>-t_-$ (resp for all $t<t+$)
\begin{equation*}
    \phi_\pm(t)=e^{-\rmi t E(V_{\pm\infty})}S(V_\pm(t))+e^{J t E(V_{\pm\infty})}e^{J t H(V_{\pm\infty})}z_{\pm\infty}+\e_\pm(t)\\
\end{equation*}
with 
\begin{eqnarray*}
\left|\dot{V}_\pm(t)+\rmi\left(E(V_\pm(t))-E(V_{\pm\infty})\right)\right|&\leq& \frac{C}{\left\langle t\right\rangle^2}\left(\|z_0\|_{H^{s}_\sigma}+\|\rho_\pm(0)\|_{H^{s}}\right)^2,\\ 
\left|V_\pm(t)-V_{\pm\infty}\right|&\leq& \frac{C}{\left\langle t\right\rangle}\left(\|z_0\|_{H^{s}_\sigma}+\|\rho_\pm(0)\|_{H^{s}}\right)^2,\\ 
\max\left\{\left\|\e_\pm(t)\right\|_{
H^{s}},\;\left\|\e_{\pm}(t) \right\|_{H^{s}_{-\sigma}},
\left\|\e_{\pm}(t) \right\|_{B^\beta_{\infty,2}}\right\}&\leq&
\ds \frac{C}{\left\langle t\right\rangle^2}\left(\|z_0\|_{H^{s}_\sigma}+\|\rho_\pm(0)\|_{H^{s}}\right)^2\\
\mbox{and }
\left\|e^{-J t H(V_{\pm\infty})}e^{J \int_0^t (E(V_{\pm}(s))-E(V_{\pm\infty})\;ds}\e_{\pm}(t) \right\|_{H^{s}_{\frac{3}{2}}}&\leq&
\ds \frac{C}{\left\langle t\right\rangle^{\frac{1}{2}}}\left(\|z_0\|_{H^{s}_\sigma}+\|\rho_\pm(0)\|_{H^{s}}\right)^2
\end{eqnarray*}
for all $t>-t_-$  (resp for all $t<t_+$).

\end{theorem}

\section{Linearized operator and exponentially stable and unstable %
manifolds}
                                                                    \label{Sec:LinearizedOp}
We study the dynamics associated with \eqref{Eq:NLD} around a
stationary state. We will use spectral properties of the linearized
operator around a stationary state.
\subsection{The spectrum of the linearized operator}
Here we study the spectrum of the linearized operator associated
with Equation~\eqref{Eq:NLD} around a stationary state~$S(U)$. Let
us recall
\begin{equation*}
H(U)=H+d^2 F(S(U))-E(U)
\end{equation*}
where~$d^2F$ is the differential of~$\nabla F$. The operator~$H(U)$
is~$\R-$linear but not $\C-$linear. Replacing~$L^2(\R^3,\C^4)$ by
$L^2(\R^3,\R^4\times\R^4)$ with the inner product obtained by taking
the real part of the inner product of~$L^2(\R^3,\C^4)$, we obtain a
symmetric operator. We then complexify this real Hilbert space and
obtain~$L^2(\R^3,\C^4\times\C^4)$ with its canonical hermitian
product. This process transforms the operator~$-\rmi$ into
\begin{equation*}
J=\left(\begin{array}{cc}
0&Id_{\C^4}\\
-Id_{\C^4}&0
\end{array}\right).
\end{equation*}
For~$\phi\in L^2(\R^3,\R^4\times\R^4)\subset
L^2(\R^3,\C^4\times\C^4)$, we still write~$\phi$ instead of
\begin{equation*}
\left(\begin{array}{c} \Re \phi\\
 \Im \phi\end{array} \right).
\end{equation*}
The extension of~$H(U)$ to~$L^2(\R^3,\C^4\times\C^4)$ is also
written~$H(U)$ and is now a real operator. The extension of $K$ (see
\eqref{Def:OperatorK}) is also written $K$. 

The linearized operator 
associated with Equation~\eqref{Eq:NLD}
around the stationary state~$S(U)$ is given by~$JH(U)$. We shall now
study its spectrum.

Differentiating~\eqref{Eq:StationaryStates}, we have that for
$U=(u_1,u_2)\in \Omega$
\begin{equation*}
{\mathcal H}_0(u_1,u_2)={\rm Span}\left\{ \frac{\p}{\p \Re
u_1}S(u_1,u_2),\; \frac{\p}{\p \Im u_1}S(u_1,u_2),\frac{\p}{\p \Re
u_2}S(u_1,u_2),\; \frac{\p}{\p \Im u_2}S(u_1,u_2)\right\}
\end{equation*}
is invariant under the action of~$JH(U)$. Differentiating  the gauge
invariance property for $S$, we notice that $JS(U)\in {\mathcal
H}_0(U)$, differentiating  the gauge invariance property for $F$, we
also obtain
\begin{equation*}
    JH(U)JS(U)=0
\end{equation*}
and differentiating \eqref{Eq:StationaryStates}, we obtain for any
$\beta\in\N^4$ with $|\beta|=1$:
\begin{equation*}
JH(U)\p_{U}^\beta S(U)=(\p_{U}^\beta E)(U)JS(U).
\end{equation*}
The space~${\mathcal H}_0(U)$ is contained in the geometric null
space of~$JH(U)$, in fact it is exactly the geometric null space as
proved in the sequel.

Now we state our results on the spectrum of $JH(U)$. The first deals
with the excited states part, we have the
\begin{proposition}\label{Prop:ExcitedSpectrumLinOp}
If Assumptions \ref{assumption:1}--\ref{assumption:Resonance} hold.
Let
\begin{equation*}
\Gamma(U)=\lim_{\substack{\e\to 0,\\\e>0}}\left\langle d^2
    F(S(U))\phi_1,\Im\left((H-\lambda_0)
    +(\lambda_1-\lambda_0)-\rmi\e\right)^{-1}P_c(H)d^2
    F(S(U))\phi_1\right\rangle
\end{equation*}
for any sufficiently small $U$. Then there exists a 
map $E_1:~B_{\C^2}(0,\e)\mapsto~\R$ with
\begin{equation*}
\begin{cases}
\Im E_1(U)=(\lambda_1-E(U))+O(|U|^4)\\
\Re E_1(U)=1/2\Gamma(U)+O(|U|^6)
\end{cases}
\end{equation*}
such that $E_1(U)$, $ \overline{E_1}(U)$, $-E_1(U)$ and
$-\overline{E_1}(U)$ are double eigenvalues of $JH(U)$ and we have $E_1(U)=E_1(|U|)$. 

For any
$s\in\R$, there exist smooth maps $k^\pm:~B_{\C^2}(0,\e)\mapsto
\left\{\left(\begin{array}{c}
\phi_1\\-i\phi_1\end{array}\right)\right\}^\bot\cap H^s$ such that
\begin{equation}\label{Identity:ExpansionExcitedStateSimple}
k_\pm(U)\\ -\left((H-E(U))+\rmi E_1(U)\right)^{-1}P_c(H)d^2
F(S(U))\frac{|U|}{\sqrt{2}}\left(\frac{u_1}{|U|}Id_{\C^4}+\frac{u_2}{|U|}K\right)\left(\begin{array}{c}
\phi_1\\-i\phi_1\end{array}\right)
\end{equation}
is $O(|U|^5)$ for any $\sigma\geq 1$.
For any $U=(u_1,u_2)\in B_{\C^2}(0,\e)$
\begin{equation*}
k_\pm(U)=\frac{|U|}{\sqrt{2}}\left(\frac{u_1}{|U|}Id_{\C^4}+\frac{u_2}{|U|}K\right)k_\pm((|U|,0))
\end{equation*} 
and defining for any $U=(u_1,u_2)\in B_{\C^2}(0,\e)$ :
\begin{equation*}
\Phi_\pm(U)=
\frac{|U|}{\sqrt{2}}\left(\frac{u_1}{|U|}Id_{\C^4}+\frac{u_2}{|U|}K\right)\left(\begin{array}{c}
\phi_1\\-i\phi_1\end{array}\right)+k_\pm(U),
\end{equation*}
we have
\begin{itemize}
    \item $\left\{\Phi_+(U),K\Phi_+(U)\right\}$ is a basis of the
    eigenspace associated with $E_1(U)$,
    \item $\left\{\overline{\Phi_+(U)},K\overline{\Phi_+(U)}
    \right\}$ is a basis of the
    eigenspace associated with $\overline{E_1(U)}$,
    \item $\left\{\Phi_-(U),K\Phi_-(U)\right\}$ is a basis of the
    eigenspace associated with $-E_1(U)$,
    \item $\left\{\overline{\Phi_-(U)},K\overline{\Phi_-(U)}\right\}$
    is a basis of the
    eigenspace associated with $-\overline{E_1(U)}$.
\end{itemize}
Moreover for any~$\beta\in \N^4$,~$s\in \R^+$
and~$p,q\in[1,\infty]$. There exist $\gamma
>0$,~$\e>0$ and~$C>0$ such that for all~$U \in B_{\C^2}(0,\e)$,
one has
\begin{equation} \label{ChapitreDeauxEq:ExpDecayExcitedState}
\|e^{\gamma \langle Q\rangle}\p_U^\beta \Phi_\pm(U)\|_{B^s_{p,q}} \leq C
\|S(U)\|_2,
\end{equation}
where $\p_{u_1,u_2}^\beta=\frac{\p^{|\beta|}}{\p^{\beta_1}\Re u_1
\p^{\beta_2}\Im u_1\p^{\beta_1}\Re u_2 \p^{\beta_2}\Im u_2}$.
\end{proposition}
\begin{proof}
For this proof, we use ideas of the proof of \cite{TsaiYau3}[Theorem
2.2]. The equation to solve for excited states is:
\begin{equation}\label{Eq:ForExcitedStateBeforeDecomposition}
    \left(JH(U)-z\right)\phi=0.
\end{equation}
Since the proof is similar for all cases, we restrict the study to $U$ 
of the form $(|U|,0)$ and (dividing by~$|U|$) to solutions the form 
$\phi=S_1+\eta$ where  $S_1$ is the normalized eigenvector of $JH$:
\begin{equation*}
S_1=\frac{1}{\sqrt{2}}\left(\begin{array}{c}
\phi_1\\-i\phi_1\end{array}\right)
\end{equation*}
and $\eta\in\left\{S_1\right\}^{{\bot}}$, the orthogonal relation is
taken in fact with respect to $J$ (but since $JS_1=\rmi S_1$, we can
take it in the usual way). For $z\in \C\setminus\rmi\R$, we obtain
the equation
\begin{equation}\label{Eq:ExcitedState}
    \eta=\left(J(H-E(U))-z\right)^{-1}P_1^{{\bot}} W(U)\left\{S_1+\eta\right\}
\end{equation}
with $P_1^{{\bot}}$ the orthogonal projector, with respect to $J$,
into $\left\{S_1\right\}^{{\bot}}$ and $W(U)=JH(U)-J(H-E(U))$. We
notice that $\left\{S_1\right\}^{{\bot}}$ is invariant under the
action of $J(H-E(U))$. To solve this equation in $\eta$ for a fixed
$u$ and $z$, we notice that if
\begin{equation*}
\Re z > 0 \mbox{ and }|\Im z| \geq m,
\end{equation*}
the series
\begin{equation*}
    k(U,z)=\left(J(H-E(U))-z\right)^{-1}P_1^{{\bot}}\sum_{k\geq 0}
    \left(-W(U)\left(J(H-E(U))-z\right)^{-1}P_1^{{\bot}}\right)^kW(U)S_1,
\end{equation*}
is convergent in $L^2$ for sufficiently small $|U|$ and $|\Im z|=O(|U|^2)$ using the Limiting
Absorption Principle~\eqref{Identity:LAP} and the bound of the resolvent $\left\|\left(H-z'\right)^{-1}\right\|\leq |\Im z'|^{-1}$ in $L^2$. Hence, we have a
solution of~\eqref{Eq:ExcitedState}.

Then we solve the equation in $z$. We obtain from Equation
\eqref{Eq:ForExcitedStateBeforeDecomposition} the equation
\begin{equation*}
    \left\langle\left( JH(U)-z\right)\phi,S_1\right\rangle=0 \mbox{ with
} \phi=S_1+k(U,z),
\end{equation*}
we infer
\begin{eqnarray*}
    z&=&\left\langle JH(U)S_1,S_1\right\rangle+\left\langle
    JH(U)k(U,z),S_1\right\rangle\\
    &=&\rmi\left(\lambda_1-\lambda_0\right)
    +\left\langle W(U)S_1,S_1\right\rangle\\
    &&+\sum_{k\geq 0}\left\langle
    JH(U)\left(J(H-E(U))-z\right)^{-1}P_1^{{\bot}}
    \left(-W(U)\left(J(H-E(U))-z\right)^{-1}
    P_1^{{\bot}}\right)^kW(U)S_1,S_1\right\rangle\\
    &=&\rmi\left(\lambda_1-\lambda_0\right)+\left\langle W(U)S_1,S_1\right\rangle\\
    &&+\sum_{k\geq 0}\left\langle
    P_1^{{\bot}}
    \left(-W(U)\left(J(H-E(U))-z\right)^{-1}
    P_1^{{\bot}}\right)^kW(U)S_1,S_1\right\rangle\\
    &&+\sum_{k\geq 0}\Bigg\langle
    \left(W(U)+z\right)\left(J(H-E(U))-z\right)^{-1}\!P_1^{{\bot}}\!\left(-W(U)\left(J(H-E(U))-z\right)^{-1}
    \!P_1^{{\bot}}\!\right)^k \!W(U)S_1,S_1\Bigg\rangle.
\end{eqnarray*}
Since $P_1^{{\bot}} S_1=0$, we introduce the function
\begin{multline*}
f(z)=\rmi\left(\lambda_1-\lambda_0\right)+\left\langle
W(U)S_1,S_1\right\rangle\\
+\sum_{k\geq 0}\left\langle
    W(U)\left(J(H-E(U))-z\right)^{-1}P_1^{{\bot}}
    \left(-W(U)\left(J(H-E(U))-z\right)^{-1}
    P_1^{{\bot}}\right)^kW(U)S_1,S_1\right\rangle.
\end{multline*}
Since  $JS_1=-\rmi S_1$, we obtain that $\Re\left\langle
W(U)S_1,S_1\right\rangle=0$, so for $z\in \C\setminus \rmi\R$ we
have
\begin{eqnarray*}
    \Re f(z)
&=& \Re \left\langle
    W(U)\left(J(H-E(U))-z\right)^{-1}P_1^{\bot}
    W(U)S_1,S_1\right\rangle
    +O(|U|^6)\\
&=& \Im \left\langle
    d^2F(S(U))\left((H-E(U))+zJ\right)^{-1}P_1^{\bot}
    d^2F(S(U))S_1,S_1\right\rangle +O(|U|^6).
\end{eqnarray*}
Then using  \eqref{FermiGoldenRule} and
\begin{equation*}
    \left((H-E(U))+zJ\right)^{-1}=\\\frac{1}{2}
    \left(\left((H-E(U))-\rmi z\right)^{-1}\left(I_{C^2}+\rmi J\right)
    +\left((H-E(U))+\rmi z\right)^{-1}\left(I_{C^2}-\rmi J\right)\right),
\end{equation*}
we obtain
\begin{multline*}
    \Im \left\langle d^2F(S(U))\left((H-E(U))+zJ\right)^{-1}
    P_1^{\bot}d^2F(S(U))S_1,S_1\right\rangle\\
    =\frac{1}{2}\Im \left\langle d^2F(S(U))\left((H-E(U))-\rmi
z\right)^{-1}d^2P_1^{\bot}F(S(U))S_1,S_1\right\rangle\\
    +\frac{1}{2}\Im \left\langle d^2F(S(U))
    \left((H-E(U))+\rmi z\right)^{-1}P_1^{\bot}d^2F(S(U))S_1,S_1\right\rangle\\
    -\Im\left\langle  d^2F(S(U))\left((H-E(U))^2+ z^2\right)^{-1}
    z JP_1^{\bot}d^2F(S(U))S_1,S_1\right\rangle,
\end{multline*}
and so using regularity results of the resolvent of
\cite{GeorgescuMantoiu}[Theorem 1.7], we obtain
\begin{multline*}
    \Im \left\langle
    d^2F(S(U))\left((H-E(U))+\left(\rmi(\lambda_1-\lambda_0)+0\right)
    J\right)^{-1}P_1^{\bot}d^2F(S(U))S_1,S_1\right\rangle\\
    =\frac{1}{2}\Im
    \left\langle  d^2F(S(U))\left((H-E(U))
    +(\lambda_1-\lambda_0)-\rmi0\right)^{-1}P_c(H)d^2F(S(U))S_1,S_1\right\rangle.
\end{multline*}
Using Assumption \ref{assumption:Resonance}, the limiting absorption
principle \eqref{Identity:LAP} and regularity results of
\cite{GeorgescuMantoiu}[Theorem 1.7], we obtain
\begin{equation*}
\Re f(z)=   1/2\Gamma(U)+O(|U|^6)
\end{equation*}
for $z$ in a ball of radius of order $|U|^2$ around
$\rmi\left(\lambda_1-\lambda_0\right)$ and for small $U$. We also
prove by the same way
\begin{equation*}
\Im f(z)=\left(\lambda_1-\lambda_0\right)+O(|U|^4)
\end{equation*}
for $z$ in a ball of radius of order $|U|^2$ around
$\rmi\left(\lambda_1-\lambda_0\right)$ and for small $U$.

So we have proved that for sufficiently small $U$, $f$ leaves a ball
around $\rmi\left(\lambda_1-\lambda_0\right)+1/2\Gamma(U)$ invariant.
With the same ideas, we prove that it is a contraction. Therefore, 
we have a fixed point $E_1(U)$ of
each $U$. Then we choose $k_+(U)=|U|k(U,E_1(U))$. Using the complex
conjugation, we obtain the eigenvalue $\overline{E_1(U)}$ and its
associated eigenspace.

The estimate on \eqref{Identity:ExpansionExcitedStateSimple} is proved usinf
$|\Re E_1(U)|=O(|U|^2)$, the Limiting
Absorption Principle~\eqref{Identity:LAP} and the bound of the resolvent 
$\left\|\left(H-z'\right)^{-1}\right\|\leq |\Im z'|^{-1}$ in $L^2$. 

Using Weyl's sequences, we prove that the essential spectrum of
$(JH(U))^*=-H(U)J$, for small $U$, is the essential spectrum of
$-HJ=-JH$. So $z$ with non zero real part is in the spectrum of
$(JH(U))^*$ if and only if it is an isolated eigenvalue. Then to
obtain $-E_1(U)$ and $-\overline{E_1(U)}$, we notice that $E_1(U)$
and $\overline{E_1(U)}$ are eigenvalues of $(JH(U))^*$. Using
the symmetry : $J(JH(U))=-(JH(U))^*J$, we show that any eigenvector
$\phi$ of $(JH(U))^*$ associated with $\lambda$, $J\phi$ is an
eigenvector of $JH(U)$ associated with $-\lambda$. Hence repeating 
the previous proof for $(JH(U))^*$, we obtain $k_-$.

The exponential decay works like in Lemma \ref{lemma:ExpDecay}
\end{proof}
\begin{remark}
If $F(z)$ is homogeneous of order $p$ then there exist
$\e,\Gamma_1,\Gamma_2>0$ such that for all $U\in B_{\C^2}(0,\e)$
\begin{equation*}
\left|U\right|^{p-2}\Gamma_1\leq\Gamma(U)\leq\left|U\right|^{p-2}\Gamma_2.
\end{equation*}
We just write $S((u_1,u_2))=u_1\phi_0+u_2K\phi_0+h((u_1,u_2))$,
expand $\Gamma(U)$ and use Assumption \ref{assumption:Resonance}
with the regularity results of the resolvente from
\cite{GeorgescuMantoiu}[Theorem 1.7].

This gives
\begin{equation*}
k_\pm(U)\\ -\left((H-E(U))+\rmi E_1(U)\right)^{-1}P_c(H)d^2
F(S(U))\frac{|U|}{\sqrt{2}}\left(\frac{u_1}{|U|}Id_{\C^4}+\frac{u_2}{|U|}K\right)\left(\begin{array}{c}
\phi_1\\-i\phi_1\end{array}\right)
\end{equation*}
is $O(|U|^{7-p})$ in $B(L^2_{-\sigma})$ for any $\sigma\in\R^+$.
\end{remark}


The following proposition deals with the essential spectrum of our
linearized operator.
\begin{proposition}\label{Prop:ContSpectrumLinOp}
If Assumptions \ref{assumption:1}--\ref{assumption:Resonance} hold.
For any sufficiently small non zero $U\in\C^2$, let
\begin{equation*}
{\mathcal H}_{1}(U)={\rm span}\left\{\Phi_+(U), K\Phi_+(U),
\overline{\Phi_+(U)}, K\overline{\Phi_+(U)}, \Phi_-(U),
K\Phi_-(U),\overline{\Phi_-(U)},K\overline{\Phi_-(U)}\right\}.
\end{equation*}
The orthogonal space  of ${\mathcal H }_0(U)\oplus {\mathcal
H}_{1}(U)$ with respect to the product associated to~$J$
\begin{equation*}
{\mathcal H}_c(U)=\left\{{\mathcal H }_0(U)\oplus {\mathcal
H}_{1}(U)\right\}^{{\bot_J}}
\end{equation*}
is invariant under the action of~$JH(U)$.

%
We also have for~${\mathbf P}_c(U)$, the orthogonal projector
onto~${\mathcal H}_c(U)$ with respect to $J$, and for~$U'\in
B_{\C^2}(U,\e)$, with sufficiently small $\e>0$, that
\begin{equation*}
\left.{\mathbf P}_c(U)\right|_{{\mathcal H}_c(U')}: {\mathcal
H}_c(U')\mapsto {\mathcal H}_c(U)
\end{equation*}
is an isomorphism and is a bounded operator
from~$H^{s}_{\sigma}(\R^3,\C^8)$ or $B^{s}_{p,q}(\R^3,\C^8)$ to
itself for any reals ~$s,\sigma\in\R$ and $p,q\in [1,\infty]$, the
inverse~$R(U',U)$ is continuous with respect to~$U$ and~$U'$ for
these norms.

Moreover, there exists $C>0$ such that we have
\begin{eqnarray}
                                                                \label{Eq:LowerBoundPc}
    \left\|\psi\right\|_X&\leq& C\left\|{P}_c(U)\psi\right\|_X,\\\;
&&\quad\forall\psi\in {\mathcal H}_c(U) \mbox{ with }
X=H^{s}_{\sigma}(\R^3,\C^8)\mbox{ or } B^{s}_{p,q}(\R^3,\C^8),\;\nonumber\\
&&\quad\forall s,\sigma\in\R,\forall p,q\in [1,\infty]\nonumber, \\
\int_{\R}\|\langle Q\rangle^{-\sigma}e^{sJH(U)}{\mathbf
P}_c(U)\psi\|^2\;ds &\leq& C\|\psi\|_2^2,\;\forall\psi\in
L^2,\;\forall \sigma\geq 1,\label{Eq:PertutbedSmoothEstimate}\\
\|e^{tJH(U)}{\mathbf P}_c(U)\psi\|&\leq& C\|\psi\|_2,\;\forall
t\in\R,\;\forall\psi\in L^2,\label{Eq:PertutbedConservationLaw}
\end{eqnarray}
and ${\mathcal H}_c(U)$ contains no eigenvector.
\end{proposition}
\begin{remark}
We use the same notation for ${\mathcal H}_c(U)$ and its real part
which appears in our main theorems. We just notice that
${\mathcal H}_c(U)$ appears when we discuss spectral properties in
our proof. Then when we talk about dynamical properties, we deal
with its real part. We remind that the real part of ${\mathcal H}_c(U)$
is left invariant by $JH(U)$.
\end{remark}
\begin{proof}
We prove that there is no other eigenvector, by proving that
smoothness estimate \eqref{Eq:PertutbedSmoothEstimate} takes place
over
\begin{equation*}
{\mathcal H}_c(U)=\left\{{\mathcal H }_0(U)\oplus {\mathcal
H}_{1}(U)\right\}^{{\bot_J}}.
\end{equation*}
First we prove that
\begin{equation*}
\left.{\mathbf P}_c((U))\right|_{{\mathcal H}_c(U')}: {\mathcal
H}_c(U')\mapsto {\mathcal H}_c(U)
\end{equation*}
is an isomorphism.
To prove it, we exhibit an inverse $R(U',U)$
which is the projector onto ${\mathcal H}_c(U')$
associated with the decomposition 
${\mathcal H}_0(U)\oplus {\mathcal H}_1(U) \oplus {\mathcal H}_c(U')$ of 
$L^2(\R^3,\C^8)$. Indeed we have
$\left\{{\mathcal H}_0(U)\oplus {\mathcal H}_1(U)\right\} \cap{\mathcal H}_c(U')=\{0\}$ 
when $U'$ and $U$ are close one to each other and
${\rm codim}{\mathcal H}_c(U')= {\rm dim}{\mathcal H}_0(U)\oplus {\mathcal H}_1(U)$, 
hence we have a decomposition of 
$L^2(\R^3,\C^8)$ into closed subspaces hence the associated
projectors are continuous. So $R(U',U)$ should be of the form
\begin{equation*}
    R(U',U)=Id+\sum_i\left|J\xi_i(U)\right\rangle
    \left\langle \alpha_i(U',U) \right|
\end{equation*}
where $\xi_i(U)$ is a basis of the eigenspaces of $JH(U)$ and
$(\alpha_i(U',U))_i$ solve the equations
\begin{equation*}
    J\xi_j(U')+\sum_i
    \left\langle J\xi_i(U),J\xi_j(U')\right\rangle\alpha_i(U',U)=0.
\end{equation*}
Such $\alpha$ exists because
the matrix $\left(\left\langle J\xi_i(U),J\xi_j(U')\right\rangle\right)_{i,j}$ 
is a Gramm matrix when $U=U'$ and
otherwise a small perturbation of such matrices for $U$ and $U'$
close one to each other and hence is invertible.

The boundedness of $R$ in ${\mathcal B}(H^{s}_{\sigma}(\R^3,\C^8))$ 
or ${\mathcal B}(B^{s}_{p,q}(\R^3,\C^8))$
follows from the exponential decay of
eigenvectors and their derivatives, the continuity of $R$ follows
from the continuity of the eigenvector with respect to the
parameters $U$ or $U'$ see Proposition
\ref{Prop:ManifoldPLS}, Lemma \ref{lemma:ExpDecay} and Proposition
\ref{Prop:ExcitedSpectrumLinOp}.

Let us now consider the orthogonal projector $P_c$ on associated
with the continuous subspace of $JH$. Since the eigenvector of $JH$
are exponentially decaying, we can extend $P_c$ to obtain an operator of
$L^2_{\pm\sigma}$ into itself. The same is true for $P_c(U)$ and hence
we can consider the extension of ${\mathcal H}_c(U)$ to
$L^2_{\pm\sigma}$. We still call it ${\mathcal H}_c(U)$. 
For all $\psi\in{\mathcal H}_c(U)$ :
\begin{equation*}
    \left\|\psi\right\|_{L^2_{-\sigma}}\leq  \left\|P_c \psi\right\|_{L^2_{-\sigma}}
+\left\|(1-P_c) \psi\right\|_{L^2_{-\sigma}}.
\end{equation*}
Since $1-P_c$ is the projector into the eigenspaces of $H$ and
$\psi$ is orthogonal to the eigenvectors of $JH(U)$, we obtain that
there exists $c>0$ such that
\begin{equation*}
 \left\|(1-P_c) \psi\right\|_{L^2_{-\sigma}}\leq c|U|
\left\|\psi\right\|_{L^2_{-\sigma}}.
\end{equation*}
Indeed, using Proposition \ref{Prop:ExcitedSpectrumLinOp}, we obtain
a $c'>0$ such that for sufficiently small non zero $U$, we have
\begin{equation*}
\left| \left\langle
\psi,J\frac{1}{\sqrt{2}}\left(\frac{u_1}{|U|}Id_{\C^4}+\frac{u_2}{|U|}K\right)\left(\begin{array}{c}
\phi_1\\-i\phi_1\end{array}\right)\right\rangle\right|\leq
\frac{1}{\left|U\right|}\left|\left\langle
\psi,J\frac{\Phi_+(U)+\Phi_-(U)}{2}\right\rangle\right| +
c'\left|U\right|\left\|\psi\right\|_{L^2_{-\sigma}}.
\end{equation*}
Hence since $\psi$ is orthogonal to $\Phi_+(U)$ and $\Phi_-(U)$, we
obtain that the projection of $\psi$ in the second eigenspace of $H$
is small, since they are invariant under the action of $J$. Using 
Proposition \ref{Prop:ManifoldPLS}, we obtain the
same thing for the first one.

Hence for a sufficiently small non zero $U$, we obtain Estimate
\eqref{Eq:LowerBoundPc} for $X=L^2_{-\sigma}$ with $\sigma>0$. The
rest of  Estimate \eqref{Eq:LowerBoundPc} follows by the same way
using the exponential decay of eigenvectors (Estimate
\eqref{ChapitreDeauxEq:ExpDecayExcitedState} and Lemma
\ref{lemma:ExpDecay}). 
We infer
\begin{eqnarray*}
\lefteqn{\|\langle Q\rangle^{-\sigma}
e^{tJH(U)}{\mathbf P}_c(U)\psi\|}\nonumber\\
&\leq& C\|\langle Q\rangle^{-\sigma}{\mathbf P}_ce^{tJH(U)}
{\mathbf P}_c(U)\psi\|\nonumber\\
&\leq& C\|\langle Q\rangle^{-\sigma}{\mathbf P}_ce^{-\rmi t(H-E(U))}
{\mathbf P}_c(U)\psi\|\nonumber\\
&&+C\|\langle Q\rangle^{-\sigma}\int_0^t {\mathbf
P}_ce^{-\rmi(t-s)(H-E(U))}D\nabla F(S(U))e^{ s JH(U)}
{\mathbf P}_c(U)\psi\,ds\|\nonumber
\end{eqnarray*}
Using estimate~\eqref{Estimate:Smoothness1} and
\eqref{Estimate:Smoothness3} of Theorem \ref{Thm:Smoothness}, we
obtain the estimate \eqref{Eq:PertutbedSmoothEstimate} for
sufficiently small~$U$:
\begin{equation*}
\int_{\R}\|\langle Q\rangle^{-\sigma}e^{-sJH(U)}{\mathbf
P}_c(U)\psi\|^2\;ds \leq C \|\psi\|^2.
\end{equation*}
Hence there is no eigenvector in the range of $P_c(U)$.
Using the inequalities \eqref{Eq:PertutbedSmoothEstimate}, the
conservation law for $H$ and Duhamel's formula :
\begin{equation*}
    e^{J t H(U)}=e^{-\rmi t
(H-E(U))}+\int_0^te^{-\rmi (t-s) (H-E(U))}Jd^2\nabla F(S(U))e^{J s
H(U)}\;ds,
\end{equation*}
we prove the estimate \eqref{Eq:PertutbedConservationLaw}.
\end{proof}

\medskip

Since ${\mathcal H}_c(U)$ is closed and ${\rm codim} {\mathcal
H}_c(U)=\dim \left\{{\mathcal H}_0(U)\oplus {\mathcal
H}_1(U)\right\}$ and $ {\mathcal H}_c(U)\cap \left\{{\mathcal
H}_0(U)\oplus\!{\mathcal H}_1(U)\right\}\!\!=~\{0\}$, we obtain $
{\mathcal H}_0(U)\oplus {\mathcal H}_1(U)\oplus{\mathcal H}_c(U)=
L^2(\R^3,\C^8)$ and the
\begin{proposition}\label{Prop:NullSpace}
Suppose that  Assumptions
\ref{assumption:1}--\ref{assumption:Resonance} hold. Then the space ${\mathcal H}_0(U)$
is the geometric null space of $JH(U)$.
\end{proposition}
\subsection{Stable, unstable and center manifold}
                                                                \label{Sec:CSUMan}
We can now obtain results 
similar to those of Bates and Jones \cite{BatesJones}. We notice 
that we won't prove that
the Cauchy problem \eqref{Eq:NLD} is locally wellposed for initial
condition outside some manifolds (built below). In fact it can be proved 
with the methods we present here or by generalizing to our case
the results of  and Vega \cite{EscobedoVega}.

We have that $JH(U)$ as an operator in $L^2(\R^3,\R^8)$ is a
closed densely defined operator that generates a continuous
semigroup on $L^2(\R^3,\R^8)$. The spectrum of $JH(U)$ in
$L^2(\R^3,\R^8)$  is the same as $JH(U)$ in $L^2(\R^3,\C^8)$ and so
it splits in three parts:
\begin{align*}
\sigma_s(U)=&\left\{\lambda\in\sigma(JH(u)),\;\Re\lambda<0\right\}
=\left\{-E_1(U),-\overline{E_1(U)}\right\}\\
\sigma_c(U)=&\left\{\lambda\in\sigma(JH(u)),\;\Re\lambda=0\right\}
=\rmi\left\{\R\setminus(-E(U),E(U))\right\}\\
\sigma_u(U)=&\left\{\lambda\in\sigma(JH(u)),\;\Re\lambda>0\right\}
=\left\{E_1(U),\overline{E_1(U)}\right\}
\end{align*}
each one is associated with a spectral real subspace, respectively
\begin{align*}
  X_s(U)&={\rm span}_\R\left\{\Re\Phi_-(U),\Im\Phi_-(U),K\Re\Phi_-(U),K\Im\Phi_-(U)\right\}\\
  X_u(U)&={\rm span}_\R\left\{\Re\Phi_+(U),\Im\Phi_+(U),K\Re\Phi_+(U),K\Im\Phi_+(U)\right\}\\
  X_c(U)&=\Re {\mathcal H}_0(U)\oplus \Re {\mathcal H}_c(U)
\end{align*}
where we used the notation $\Re
\Psi=(1/2)\left(\Psi+\overline{\Psi}\right)$ and
$\Im\Psi=-(\rmi/2)\left(\Psi-\overline{\Psi}\right)$ and $\Re
X=\left\{\Re \Psi,\Psi\in X\right\}$, the real part of the space
$X$. The spaces $X_s(U)$ and $X_u(U)$ are finite dimensional. Let us
write $\pi^c(U)$, $\pi^s(U)$ and $\pi^u(U)$ for the projector
associated with the decomposition $X_c(U)\oplus X_s(U)\oplus
X_u(U)$. Since the eigenvectors belongs also to $L^2_\sigma$ for any
$\sigma\in \R$, the projector $P_c(U)$ and $\pi^c(U)$ can be defined
in $L^2_\sigma$ for any real $\sigma$. We can extend, by this way,
the spaces ${\mathcal H}_c(U)$ and $X_c(U)$ to $L^2_\sigma$ for any
$\sigma\in \R$. We have the
\begin{lemma}
If Assumptions \ref{assumption:1}--\ref{assumption:Resonance} hold.
Then any $\sigma\in \R$, there exist $r,C_1,C_2>0$ such that for all
$t\in \R$, we have
\begin{align}
  C_1  e^{-\gamma(U)t}\leq \left\|e^{tJH(U)}\pi^s(U)\right\|_{{\cal B}(L^2_\sigma)}&\leq C_2
  e^{-\gamma(U)
  t},\label{Estimate:NormLinearizedFlow1}\\
C_1 e^{\gamma(U)t}\leq \left\|e^{tJH(U)}\pi^u(U)\right\|_{{\cal
B}(L^2_\sigma)}&\leq C_2
e^{\gamma(U) t},\label{Estimate:NormLinearizedFlow2}\\
\left\|e^{tJH(U)}\pi^c(U)\right\|_{{\cal B}(L^2_\sigma)}&\leq C_2
  \langle t\rangle^{r},\label{Estimate:NormLinearizedFlow3}
\end{align}
where $\gamma(U)=\Re E_1(U)$.
\end{lemma}
\begin{proof}
The statements for the spaces $X^s(U)$ and $X^u(U)$ follows
from \eqref{ChapitreDeauxEq:ExpDecayExcitedState}.

The statement about $X^c(U)$ is a little more complicate. We notice
that we are not
looking for an optimal $r$. 

First, the result for $e^{-\rmi t (D_m+V)}$ in $L^2_\sigma$
with $\sigma\in 2\N$ follows from \cite[Theorem 8.5]{Thaller} 
(see also Proposition \ref{Prop:WeightedEstimates} below),
which is based on the charge conservation. The case $\sigma\in \R$ follows by duality and
interpolation. 

Then for $e^{tJH(U)}\pi^c(U)$, we use Duhamel's formula :
\begin{equation*}
    e^{J t H(U)}\pi^c(U)=e^{-\rmi t
(H-E(U))}\pi^c(U)+\int_0^te^{-\rmi (t-s) (H-E(U))}Jd^2\nabla
F(S(U))e^{J s H(U)}\pi^c(U)\;ds,
\end{equation*}
then the assertion for $e^{tJH(U)}\pi^c(U)$ follows from the
assertion for $e^{-\rmi t (D_m+V)}$, the charge conservation of
$e^{tJH(U)}P_c(U)$ (see \eqref{Eq:PertutbedConservationLaw}), the
fact that $e^{tJH(U)}S(U)=S(U)$, $e^{tJH(U)}\p_U^\beta
S(U)=\p_U^\beta S(U)+t\p_U^\beta E(U)S(U)$ and Lemma
\ref{lemma:ExpDecay}.
\end{proof}

{\it By now we do not restrict our study to the space
$L^2(\R^3,\R^8)$, we extend it to $L^2_\sigma(\R^3,\R^8)$ for any
$\sigma\in \R$, but we still write ${\mathcal H}_c(U)$ and $X_c(U)$
for the extensions of these spaces to $L^2_\sigma(\R^3,\R^8)$ for
any $\sigma\in \R$.}

We now study the behavior of the solutions in $L^2_\sigma$ of
\eqref{Eq:NLD} centered around $S(U)$:
\begin{equation}\label{Eq:Centered}
    \p_t\phi=JH(U)\phi +JN(U,\phi)
\end{equation}
where $H(U)=H+d^2F(S(U))-E(U)$ and $N(U,\phi)=\nabla
F(S(U)+\phi)-\nabla F(S(U))-d^2 F(S(U))\phi$ and $d^2F$ is the
differential of~$\nabla F$.

In this subsection, we study a modified equation which coincides
with \eqref{Eq:Centered} as long as the solution stays in a neigborhood of a small $S(U)$:
\begin{equation}\label{Eq:CenteredMod}
    \p_t\phi=JH(U)\phi +JN_\e(U,\phi)
\end{equation} 
where $N_\e(U,\eta)=\rho(\e^{-1}\eta)N(U,\eta)$ and $\rho$ is a
smooth function with compact support around $0$.


We state the
\begin{proposition}[Center-Stable Manifold]                   \label{Prop:CentrStabMan}
If Assumptions \ref{assumption:1}--\ref{assumption:Resonance} hold.
Then for any sufficiently small non zero $U$, there exists around
$S(U)$ a unique invariant smooth center-stable manifold $W^{cs}(U)$
for \eqref{Eq:CenteredMod} build as a graph with value in $X_u(U)$
and tangent to $S(U)+X_c(U)\oplus X_s(U)$ at $S(U)$.

Any solution $\phi\in L^2_\sigma$ of \eqref{Eq:CenteredMod}
initially in the neighborhood of $S(U)$ tends as $t\to -\infty$ to
$W^{cs}(U)$ with
\begin{equation*}
dist_{L^2_\sigma}(\phi(t), W^{cs}(U))=O(e^{\gamma t}) \mbox{ as }t
\to -\infty
\end{equation*}
for any $\gamma \in (0,\gamma(U))$, any $s,\sigma\in\R$ and for any
sufficiently small neighborhood $V$ of $S(U)$ any solution in $V$
not in $W^{cs}(U)$ leaves $V$ in finite positive time.
\end{proposition}
\begin{remark}                                                  \label{Rem:TauGamma}
For any $s\in \R^+$, due to the exponential decay of eigenvectors,
even if $\phi\notin H^s_\sigma$, there exists $\psi\in W^{cs}(U)$
such that $\phi-\psi\in H^s_\sigma$ and we have
\begin{equation*}
dist_{H^s_\sigma}(\phi(t), W^{cs}(U))=O(e^{\gamma t}) \mbox{ as }t
\to -\infty
\end{equation*}
as shown in the following proof.

If we only consider small solutions, we obtain a locally invariant
manifold for the equation \eqref{Eq:Centered}, that is to say that
for any initial condition in the manifold there exist a
corresponding solution of \eqref{Eq:Centered} which stays in this
manifold in a small interval of time around $0$. We notice that in
the following proofs the size of this invariant manifold, which is
given by $\e$, is a function of $U$ and this function is
$O(\gamma(U))$. By now, we call this function $r$.
\end{remark}
\begin{proof}
Our proof is an adaptation of the one of Bressan \cite{Bressan} and
we refer to it for more details. We make the proof only for the case
$\sigma=0$, the proof in the general case is similar.

First we prove that there is a global solution of the equation
\eqref{Eq:CenteredMod} which do not grow much as $t\to +\infty$. We
look for solution as a fixed point:
\begin{equation*}
    y(t)={\mathcal G}_\e(y_0,y)(t)
\end{equation*}
for any $y_0\in X_s(U)\oplus X_c(U)$ where for small positive $\e$
\begin{multline*}
    {\mathcal G}_\e(y_0,\eta)(t)=e^{tJH(U)}y_0+\int_0^t
    e^{(t-s)JH(U)}\pi^c(U)JN_\e(U,\eta(s))\;ds\\
    +\int_0^t
    e^{(t-s)JH(U)}\pi^s(U)JN_\e(U,\eta(s))\;ds-\int_t^{+\infty}
    e^{(t-s)JH(U)}\pi^{u}(U)JN_\e(U,\eta(s))\;ds,
\end{multline*}
with $\pi^*(U)$ the projector into $X^*(U)$ with
respect to the decomposition
$\ds\oplus_{*\in\{c,s,u\}}X^{*}(U)$.

Let us introduce
for $\gamma(U)=\Re E_1(U)$ and any $\Gamma$ smaller than
$\gamma(U)$, the space
\begin{equation*}
Y_\Gamma=\left\{y:\R\mapsto L^2(\R^3,\C^4),\;\exists C>0,
\;\|y(t)\|_2\leq C e^{\Gamma|t|},\;\forall t\in\R\right\}.
\end{equation*}
For sufficiently small $\e>0$, the map ${\mathcal G}_\e(y_0,\cdot)$
leaves $Y_\Gamma$ invariant and is continuous for the norm
\begin{equation*}
N_\Gamma : y\mapsto \sup_{t\in\R}\left\{\|y(t)\|_2e^{-\Gamma|t|}\right\}.
\end{equation*}
Moreover, it is a strict contraction for sufficiently small $U$ and
$\e>0$. Actually we choose $\e$ as a function of $\Gamma$ which is 
$O(\Gamma)$. In fact since $Y_\Gamma\subset Y_{\Gamma'}$ for $\Gamma<
\Gamma'$, we obtain that $\e$ as a function of $U$ is a
$O(\gamma(U))$. This proves the existence of the fixed point $y$.

Then we fix $h^{cs}_{U}(y_0)=y(0)-y_0$. The invariance of the graph of
$h^{cs}_{U}$ by the flow of Equation \eqref{Eq:CenteredMod} is
immediate.

Now we prove the smoothness property. We have $N_\e(U,\eta)$ is $l$
times differentiable in $\eta$ from $Y_{\Gamma'}$ to $Y_\Gamma$ if
$(l+1)\Gamma'\leq \Gamma$ and ${\mathcal G}_\e$ is $l$ times
differentiable from $X_c(U)\oplus X_s(U)\times Y_{\Gamma''}$ to
$Y_\Gamma$ if $2l\Gamma''\leq \Gamma$ (see \cite{Bressan}). We
introduce the family $(\eta_n)_{n\in\N}$ satisfying
\begin{equation*}
    \eta_0 = 0 \mbox{ and }
    \eta_{n+1} = {\mathcal G}_\e(y_0,\eta_n).
\end{equation*}
This sequence converge to $y$ (the fixed point) in $Y_\Gamma$.
Moreover, as functions of $y_0$, the convergence is uniform in
$Y_{\Gamma}$ (endowed with the norm $N_\Gamma$) on bounded sets 
of $X_s(U)\oplus X_c(U)$.

We want to prove that the sequence of their derivatives of order $k$
with respect to $\eta$ also converges in $Y_{\Gamma}$ on bounded
sets for any $\Gamma < \gamma(U)$. We prove it by induction in $k$.
So suppose that $(\p^j \eta_n)_{n\in\N}$ is converging in
$Y_{\Gamma}$ for all $j< k$ and any $\Gamma < \gamma(U)$. Then we
have that (see \cite{Bressan})
\begin{align*}
    \p \eta_n&=\p  {\mathcal G}_\e(y_0,\eta_{n-1})=L+
    M\left(\p_\eta N_\e(U,\eta_{n-1})\p \eta_{n_1}\right)\\
    \p^k \eta_n&=\p^k  {\mathcal G}_\e(y_0,\eta_{n-1})=
    M\left(\p_\eta N_\e(U,\eta_{n-1})\p^k
    \eta_{n-1}+\Psi_k(\eta_{n-1},\ldots,\p^{k-1}\eta_{n-1})\right),\;\forall k\geq 2
\end{align*}
with $L=e^{tJH(U)}$ and
\begin{equation*}
(M\eta)(t)=-\int_0^te^{(t-s)JH(U)}\pi^{cs}(U)J\eta(s)\;ds
+\int_t^{+\infty}e^{(t-s)JH(U)}\pi^{u}(U)J\eta(s)\;ds
\end{equation*}
and $\Psi_k$ a smooth function of $k$ parameter. Hence since
$M\circ\p_\eta N_\e(U,y_{n-1})$ is a strict contraction in
$Y_{\Gamma}$ for sufficiently small $\e$ and $U$ (once more $\e$ is
a $O(\gamma(U))$), this proves the convergence of the sequence of
$k$-th derivatives in $Y_{\Gamma}$ on bounded sets for any $\Gamma <
\gamma(U)$. Hence the sequences of derivatives of
$(\eta_n)_{n\in\N}$ in $Y_{\Gamma}$ on bounded sets for any $\Gamma
< \gamma(U)$. This gives the differentiability at any order of
$y(0)=h(y_0)$. This also gives, since $N(U,\eta)=O(|\eta|^2)$ around
zero, that $h(y_0)=O(|y_0|^2)$ around zero.

Now we want prove that $W^{cs}(U)$ is attractive in negative time.
In fact $W^{cs}(U)$ is the graph of a smooth function
$h:X^{cs}\mapsto X^u(U)$. Let $\eta$ be such that $S(U)+\eta$ is a
solution of \eqref{Eq:NLD}, we have
\begin{equation*}
\p_t \eta=JH(U)\eta+JN_\e(U,\eta).
\end{equation*}
\begin{equation*}
    \eta=y+r=y+h(y)+z
\end{equation*}
with $y=\pi^{cs}(U)\eta$ and we have the following equation for
$z\in X^u(U)$
\begin{equation*}
\p_t z= JH(U)z + M(U,y,z)
\end{equation*}
where
\begin{equation*}
M(U,y,z)=\pi^u(U)\left\{JN_\e(U,\eta)-JN_\e(U,y+h(y))\right\}\\
-Dh(y)\pi^{cs}(U)\left\{JN_\e(U,\eta)-JN_\e(U,y+h(y))\right\}.
\end{equation*}
Using Duhamel's formula, we obtain
\begin{equation*}
    z(t)=e^{tJH(U)}z(0)+\int_0^te^{(t-s)JH(U)}M(U,y(s),z(s))\;ds.
\end{equation*}
We obtain since $z\in X^u(U)$
\begin{equation*}
    \left\|z(t)\right\|\leq e^{\gamma(U)
    t}\left\|z(0)\right\|+C\int_0^te^{(t-s)\gamma(U)}\left\|M(U,y(s),z(s))\right\|\;ds
\end{equation*}
and so for $\gamma \in (0,\gamma(U))$
\begin{equation*}
e^{-\gamma t}\left\|z(t)\right\|\leq
\left\|z(0)\right\|+C|U|\sup_{s\in[0,t]}\left\{e^{-\gamma
s}\left\|z(s)\right\|\right\}+o(\sup_{s\in[0,t]}\left\{e^{-\gamma
s}\left\|z(s)\right\|\right\})
\end{equation*}
where $C$ do not depend of $U$ and $z$. Hence if $z(0)$ and $U$ are
small, we have that there exists $c>0$ such that
$\left\|z(t)\right\|\leq c e^{\gamma t}$ for all $t\leq 0$. We
notice that since $X^u(U)\subset H^s_\sigma$ for any $s\in\R^+$ and
is finite dimensional (see Lemma \ref{Prop:ExcitedSpectrumLinOp}),
the time decay in $L^2_\sigma$ gives also a time decay in
$H^s_\sigma$ for any $s\in\R^+$.

Now choose $V$ a sufficiently small neighborhood of $0$ and $\phi$ a
solution of \eqref{Eq:CenteredMod} initially in $V$ but not in
$W^{cs}(U)$. Suppose it stays in $V$ in positive time. We obtain
that $\phi\in Y_\Gamma$. We have
\begin{multline*}
    \phi(t)=e^{tJH(U)}\left(\pi^s(U)+\pi^c(U)\right)\phi(0)\\+\int_0^t
    e^{-sJH(U)}\pi^s(U)JN_\e(U,\phi(s))\;ds
    +\int_0^t
    e^{-sJH(U)}\pi^c(U)JN_\e(U,\phi(s))\;ds\\
    +e^{tJH(U)}\pi^u(U)\left(\pi^u(U)\phi(0)+\int_0^{\infty}
    e^{-sJH(U)}\pi^u(U)JN_\e(U,\phi(s))\;ds\right)\\
    -
    \int_t^{\infty}
    e^{-sJH(U)}\pi^{u}(U)JN_\e(U,\phi(s))\;ds.
\end{multline*}
with
\begin{equation*}
\pi^u(U)\phi(0)+\int_0^{\infty}
e^{-sJH(U)}\pi^u(U)JN_\e(U,\phi(s))\;ds\neq 0.
\end{equation*}
Hence we obtain with \eqref{Estimate:NormLinearizedFlow2}, that
$\phi(t)$ exponentially tends to infinity in norm. This is a
contradiction so $\phi$ leaves $V$ in finite time.
\end{proof}

Then reversing the time direction that is to say replacing $H$ by
$-H$ and $F$ by $-F$, we obtain with this theorem a locally
invariant center unstable manifold with the corresponding
properties:
\begin{proposition}[Center-Unstable Manifold]                 \label{Prop:CentrUnstabMan}
If Assumptions \ref{assumption:1}--\ref{assumption:Resonance} hold.
Then for any sufficiently small non zero $U$, there exists around
$S(U)$ a unique smooth invariant center unstable manifold
$W^{cu}(U)$ for \eqref{Eq:CenteredMod}, build as a graph with value
in $X_S(U)$ and tangent to $S(U)+X_c(U)\oplus X_u(U)$ at $S(U)$.

%
Any solution $\phi\in L^2_\sigma$ of \eqref{Eq:CenteredMod}
initially in the neighborhood of $S(U)$ tends as $t\to +\infty$ to
$W^{cu}(U)$ with for any $s\in\R^+$
\begin{equation*}
dist_{H^s_\sigma}(\phi(t), W^{cu}(U))=O(e^{-\gamma t}) \mbox{ as }t
\to +\infty
\end{equation*}
and for $\gamma \in (0,\gamma(U))$, any $s,\sigma\in\R$ and for any
$V$ sufficiently small neighborhood of $S(U)$ any solution in $V$
not in $W^{cu}(U)$ leaves $V$ in finite negative time.

\end{proposition}

We can build by the same way a center manifold which is the
intersection of the previous:
\begin{proposition}[Center Manifold]                          \label{Prop:CentrMan}
If Assumptions \ref{assumption:1}--\ref{assumption:Resonance} hold.
Then for any sufficiently small non zero $U$, there exists around
$S(U)$ a unique smooth invariant center manifold $W^{c}(U)$ for
\eqref{Eq:CenteredMod}, build as a graph with value in $X_s(U)\oplus
X_u(U)$ and tangent to $S(U)+ X_c(U)$ at $S(U)$.

Moreover, we have $W^c(U)=W^{cs}(U)\cap W^{cu}(U)$ and $W^c(U)$
contains the part of the PLS manifold which is in a small
neighborhood of $S(U)$.
\end{proposition}
\begin{proof}
We build the center manifold with the same method as in the previous
cases. We can also build a center-unstable manifold inside
center-stable manifold. More precisely, let $h^{s}_{U}: X_c(U)\oplus
X_s(U)\mapsto X_u(U)$ be the map defining center-stable manifold and
$h^{u}_{U}: X_c(U)\oplus X_u(U)\mapsto X_s(U)$ be the map defining
center-unstable manifold. A solution $y=S(U)+y_c+y_s+y_u$ with
$y_*\in X_*(U)$ for $*\in\{c,s,u\}$ is in the center-stable manifold
if $y_u=h_{U}^s(y_c,y_s)$. Hence to obtain a center-unstable manifold
inside center-stable manifold one has to solve, for each $y_c$, the
equation
\begin{equation*}
    y_s=h_{U}^u(y_c,h_{U}^s(y_c,y_s)),
\end{equation*}
which can be solve inside a small ball for small $y_c$ and small $U$
by means of the fixed point theorem, since $h_{U}^*(y_c,z)$ is a
$O(|y_c|^2+|z|^2)$ around zero for $*\in\{s,u\}$.

By the same way, we can also build a center-stable manifold inside
the center-unstable manifold.

Using the uniqueness of the center manifold, we obtain that this two
manifolds are equal to the center manifold and  $W^c(U)=
W^{cs}(U)\cap W^{cu}(U)$.

Then any stationary states in a small neighborhood of $S(U)$
converges to $W^{cs}(U)$ and $W^{cu}(U)$ using the stabilization
results of the Proposition \ref{Prop:CentrStabMan} and Proposition
\ref{Prop:CentrUnstabMan}. Hence, we have that it belongs to
$W^{cs}(U)\cap W^{cu}(U)=W^c(U)$.
\end{proof}

In the two following sections, we study the dynamic respectively
inside and outside the center manifold.

\section{The dynamic inside the center manifold}
                                                                \label{Sec:Stabilization}
In this section, we prove that the dynamic inside the center
manifold around $S(V_0)$, for small non zero $V_0$, relaxes towards
the PLS manifold. To this end, we use Theorem \ref{Thm:Smoothness}
and Theorem \ref{Thm:Strichartz} about the time decay of the
propagator associated with $H$.

\subsection{Decomposition of the system}

Like in \cite{NBStableDirectionSmallSolitonNR}, we decompose a
solution~$\phi\in W^c(V_0)$ of the equation \eqref{Eq:NLD} with
respect to the spectrum of~$JH(U)$, with $U$ specified in the
sequel, and we study the equations for these different parts of the
decomposition. We introduce
\begin{multline*} {\mathcal
H}_0^{\bot_J}(u_1,u_2)=\Bigg\{\eta\in L^2(\R^3,\C^8), \left\langle
J\eta, \frac{\p}{\p \Re u_1}S(u_1,u_2)\right\rangle=\left\langle
J\eta, \frac{\p}{\p \Im
u_1}S(u_1,u_2)\right\rangle\\
=\left\langle J\eta, \frac{\p}{\p \Re
u_2}S(u_1,u_2)\right\rangle=\left\langle J\eta, \frac{\p}{\p \Im
u_2}S(u_1,u_2)\right\rangle=0\Bigg\}.
\end{multline*}
In fact, we have
\begin{equation*}
{\mathcal H}_0^{\bot_J}(U)={\mathcal H}_{1}(U)\oplus {\mathcal
H}_c(U)
\end{equation*}
which is invariant under the action of~$JH(U)$. We recall that
${\mathcal H}_{1}(U)$ is defined in Proposition
\ref{Prop:ExcitedSpectrumLinOp} and ${\mathcal H}_c(U)$ in
Proposition \ref{Prop:ContSpectrumLinOp}. We have the
\begin{lemma}
                                                                \label{Lem:DecompositionLemma}
If Assumptions \ref{assumption:1}--\ref{assumption:NonLinearity}
hold. Let $s,\sigma\in\R$ there exist $\e,\e'>0$ such that for the
manifold
\begin{equation*}
    \Sigma=\left\{(U,\eta),\;U
\in B_{\C^2}(0,\e'),\;\eta \in {\mathcal H}_0^{\bot_J}(U)\right\}
\end{equation*}
endowed with the metric of $\C^2\times H^{s}_{\sigma}$ and any
function $\phi\in B_{H^{s}_{\sigma}}(0,\e)$, there exist a unique
$(U,\eta)\in\Sigma$ with
\begin{equation*}
\phi=S(U) + \eta.
\end{equation*}
Moreover, there exists a neighborhood ${\mathcal O}$ of $(0,0)\in
\Sigma$ such that the mapping $\phi\mapsto (U, \eta)\in {\mathcal
O}$ is smooth.
\end{lemma}
\begin{proof}
The prove that $\Sigma$ is manifold, we use Proposition
\ref{Prop:ContSpectrumLinOp} which gives that it is locally isomorph
to some open subset of $\C^2\times {\mathcal H}_c$ endowed with the
metric of $\C^2\times H^{s}_{\sigma}$. Then this is a consequence of
the inverse function theorem like in \cite[Lemma
2.3]{GustafsonNakanishiTsai}.
\end{proof}
 For any solution $\phi$
of~\eqref{Eq:NLD} on an interval of time~$I$ containing $0$, we
write for~$t\in I$
\begin{equation*}
\phi(t)=e^{-\rmi\int_0^tE(U(s))\,ds}\left(S(U(t)) + \eta(t)\right).
\end{equation*}
where~$\eta(t)\in {\mathcal H}_0^{\bot_J}(U(t))$ and we want to solve the
equation
\begin{equation}
                                                                \label{Eq:FirstForEta}
\begin{array}{ccl}
\rmi\p_t \eta&=&\left\{H-E(U)\right\}\eta+ \left\{\nabla
F(S(U)+\eta)-\nabla
F(S(U))\right\}-\rmi dS(U)\dot{U}\\
&=&\left\{H+d^2F(S(U))-E(U)\right\}\eta+N(U,\eta)-\rmi dS(U)\dot{U}
\end{array}
\end{equation}
for~$\eta(t)\in {\mathcal H}_0^{\bot_J}(U(t))$. Here~$d^2F$ is the
differential of~$\nabla F$ and~$dS$ the differential of~$S$ in
$\R^2$. To close the system, we need the equation for~$U$. This follows
from the condition
\begin{equation*}
\langle \eta(t),JdS(U(t))\rangle=0.
\end{equation*}
After a time derivation (like in\cite{NBStableDirectionSmallSolitonNR}),
we obtain the equation:
\begin{equation*}
\dot{U}(t)=-A(U(t),\eta(t))\langle N(U(t),\eta(t)) ,dS(U(t))\rangle.
\end{equation*}
where
\begin{equation*}
A(U,\eta)=[\langle JdS(U),dS(U)\rangle - \langle
J\eta,d^2S(U)\rangle]^{-1}
\end{equation*}
the matrix ~$[\langle JdS(U(t)),dS(U(t))\rangle - \langle
J\eta(t),d^2S(U(t))\rangle]$ is invertible for small~$|U(t)|$ and
$\|\eta(t)\|_2$ since we have
\begin{equation*}
[\langle JdS(U(t)),dS(U(t))\rangle - \langle
J\eta(t),d^2S(U(t))\rangle] = \left(
\begin{array}{cc}
J&0_2\\
0_2&J
\end{array}\right) + O(\left|U(t)\right|+\|\eta(t)\|_2).
\end{equation*}

\begin{lemma}
                                                                \label{Lem:DecompositionCentralManifofd2}
For any $s,s',\sigma\in\R$,  any $p,q\in[1,\infty]$, any
$V_0\in\C^2\setminus \{0\}$ sufficiently small there exist
$\e,\e'>0$, such that for the manifold
\begin{equation*}
    {\mathcal S}(V_0,\e)=\left\{(U,z);\;U\in B_{\C^2}(V_0,\e),z\in{\mathcal H}_c(U)\cap
B_{H^{s'}_\sigma}(0,\e')\right\},
\end{equation*}
endowed with the metric of $\C^2\times H^{s'}_{\sigma}$, there
exists a unique map $g: S(V_0,\e)\mapsto B^s_{p,q}(\R^3,\C^4)$ which
is smooth and satisfies for all $(U,z)\in S(V_0,\e)$, $g(U,z)\in
{\mathcal H}_1(U)$, $z+g(U,z)\in{\mathcal H}_0^{\bot_J}(U)$, and
$S(U)+z+g(U,z)\in W^c(V_0)$. Moreover, we have
$\left\|g(U,z)\right\|_{B^s_{p,q}}=O(\left\|z\right\|_{H^{s'}}^2)$.
\end{lemma}
\begin{proof}
The fact that ${\mathcal S}(V_0,\e)$ is  manifold here is proved
like in Lemma \ref{Lem:DecompositionLemma}.

Then if $h^c$ is the function for which $W^c(V_0)$ is the graph. Any
$\phi\in L^2(\R^3,\R^8)$ can be written in the form
$S(V_0)+\tilde{U}\cdot DS(V_0)+\xi+\rho$ with $\rho\in {\mathcal
H}_1(V_0)$ and $\xi\in {\mathcal H}_c(V_0)$. It can be also written
in the form $S(U)+z+r$ with $r \in{\mathcal H}_1(U)$ and
$z\in{\mathcal H}_c(U)$. These two decompositions in fact defines
two bijective smooth maps in sufficiently small sets (for the first
we have a linear decomposition, for the second see Lemma
\ref{Lem:DecompositionLemma}). We write $\Psi$ for the first and
$\Phi$ for the second. Then $f=\Psi\circ \Phi^{-1}$ has $3$
components following the decomposition ${\mathcal
H}_0(V_0)\oplus{\mathcal H}_1(V_0)\oplus{\mathcal H}_c(V_0)$, we
write them $(f_1,f_2,f_3)$. Then $g$ is the solution of the implicit
equation in $r$
\begin{equation*}
    F(U,z,r)=f_2(U,z,r)-h^c(f_1(U,z,r),f_3(U,z,r))=0
\end{equation*}
which can be solved by the implicit function theorem in
$H^{s'}_{\sigma}$ since $\p_r F(V_0,0,0)$ is invertible from
${\mathcal H}_1(V_0)$ to itself because $\p_r f_2(V_0,0,0)$
($f_2(V_0,r,0)=r$) is invertible from ${\mathcal H}_1(V_0)$ to
itself and $D h_c(0,0)=0$.
%

The  smoothness of $g$ in the Besov spaces follows from the fact
that $g(U,z)\in {\mathcal H}_1(U)$  and the exponential decay for
excited states and their derivatives given by
\eqref{ChapitreDeauxEq:ExpDecayExcitedState}.

Then we notice that for any $U$ close to $V_0$, the previous proof
can be applied to $W^c(U)$. It shows that $W^c(U)$, $W^c(V_0)$ are
in a neighborhood of $S(V_0)$ the graph of a two functions on
${\mathcal S}(V_0,\e)$ equal up to a translation in $\C^2$ of the
first argument. Hence their graphs are equal, so locally
$W^c(U)=W^c(V_0)$. The last assertion then follows from the fact
that at $S(U)$, $W^c(U)$ is tangent to $S(U)+X^c(U)$ and $X^c(U)$ is
orthogonal to ${\mathcal H}_1(U)$.
\end{proof}

Hence decomposing~$\eta$ with respect to the spectrum of
$JH(U)$,
we write
\begin{equation*}
\eta(t)=g(U(t),z(t))+z(t)
\end{equation*}
with~$z\in {\mathcal H}_c(U)\cap L^2(\R^3,\R^8)$.
We obtain the system
\begin{equation*}
\begin{cases}
\dot{U}=-A(U,\eta)\langle N(U,\eta) ,dS(U)\rangle \\
\p_t z=JH(U)z+{\mathbf P}_c(U)JN(U,\eta)\\
\qquad+{\mathbf P}_c(U(v))dS(U(v))A(U(v),\eta(v)) \langle
N(U(v),\eta(v)) ,dS(U(v))\rangle+(dP_c(U))A(U,\eta)\langle N(U,\eta)
,dS(U)\rangle\eta
\end{cases}
\end{equation*}
with
\begin{equation*}
\eta(t)=z(t) +g\left(U(t),z(t)\right).
\end{equation*}
We notice that this equation is defined only \underline{for~$z$
small with real values and~$U$ small}. We now study this system.

\subsection{The stabilization towards the PLS manifold}
                                                                \label{Section:Stabilization}

We now show that any solution of \eqref{Eq:NLD} which belongs to the
center manifold $W^c(V_0)$, for a small non zero $V_0$, stabilizes
as $t\to \pm\infty$ towards the manifold of the stationary states
inside $W^c(V_0)$. To this end, we will use Theorem
\ref{Thm:Smoothness} and Theorem~\ref{Thm:Strichartz} to prove that
$z$ tends to zero in some sense.

Let us define for any~$\e,\,\delta>0$
\begin{equation*}
{\mathcal U}(\e,\delta)\\
=\left\{U\in{\mathcal
C}^{1}(\R,B_{\C^2}(V_0,\e)),\;\|\dot{U}\|_{L^1(\R)\cap
L^\infty(\R)}\leq\delta^2\right\}
\end{equation*}
and for any~$U\in {\mathcal U}(\e,\delta)$, let~\underline{$s,\beta$
be such that~$s'> \beta+2>2$ and~$\sigma>3/2$},
\begin{multline*}
{\mathcal Z}(U,\delta)=\Bigg\{z\in {\mathcal C}(\R,
L^2(\R^3,\R^8)),\;z(t)\in {\mathcal
H}_c(U(t)),\,\\\max\left[\|z\|_{L^\infty(\R,H^{s})},
\|z\|_{L^2(\R,H^{s}_{-\sigma})},
\|z\|_{L^2(\R,B^\beta_{\infty,2})}\right]\leq\delta\Bigg\},
\end{multline*}
and $\e,\delta$ are small enough to ensure that for $U\in{\mathcal
U}(\e,\delta)$ and $z\in {\mathcal Z}(U,\delta)$
\begin{equation*}
    \ds S(U)+z+g(U,z)\in W^c(V_0)\cap B_{H^{s}}(S(V_0),r(V_0)),
\end{equation*}
where $g$ is defined by Lemma
\ref{Lem:DecompositionCentralManifofd2} and $r$ in Remark
\ref{Rem:TauGamma}. It will appear later that $\delta$ is of the 
same order as $\left\|z_0\right\|_{H^{s}}$ (see Lemma 
\ref{Lem:StabilizationForT} below).

\subsubsection{some useful lemma}

In the rest of our study, we will need some technical lemmas, which we collect here.
\begin{lemma}
                                                                \label{Lem:EstimateOnN}
If Assumptions \ref{assumption:1}--\ref{assumption:NonLinearity} hold.
Let~$\sigma,\sigma'\in\R$,~$s>1$
and~$p,\widetilde{p_1},p_1,p_2,q\in[1,\infty]$ such that
\begin{equation*}
\frac{1}{p}+\frac{s}{3}\geq\frac{1}{p_1}+\frac{1}{p_2}\geq
\frac{1}{p}.
\end{equation*}
and
\begin{equation*}
\frac{1}{p}+\frac{s}{3}\geq\frac{1}{\widetilde{p_1}}.
\end{equation*}
Then there exist~$\e>0$ and~$C>0$ such that for all~$U\in
B_\C(0,\e)$ and~$\eta \in B_{p_2,q}^s(\R^3,\R^8)\cap
L^\infty(\R^3,\R^8)$ with $\langle Q\rangle^\sigma \eta \in
B_{p_1,q}^s(\R^3,\R^8)$ and $\langle Q\rangle^{\sigma'} \eta \in
B_{\widetilde{p_1},q}^s(\R^3,\R^8)$, we have
\begin{multline}
                                                                \label{Estimate:OnN2-1}
\left\|\langle Q\rangle^\sigma N(U,\eta)\right\|_{B_{p,q}^s} \leq
C\left(s,F,\left|U\right|+\left\|\eta\right\|_{L^\infty} \right)
\left|U\right| \left\|\eta\right\|_{L^\infty} \left\|\langle
Q\rangle^{\sigma'}\eta\right\|_{B_{\widetilde{p_1},q}^s}\\
+C\left(s,F,\left|U\right|+\left\|\eta\right\|_{L^\infty\cap
B_{p_2,q}^s} \right) \left\|\eta\right\|_{L^\infty}^2 \left\|\langle
Q\rangle^\sigma\eta\right\|_{B_{p_1,q}^s}.
\end{multline}
\end{lemma}
\begin{proof}
We recall the definition
\begin{equation*}
N(U,\eta)=\nabla F(S(U)+\eta)-\nabla F(S(U))-d^2F(S(U)) \eta.
\end{equation*}
We have
\begin{equation*}
N(U,\eta)=\int_0^1\int_0^1
d^3F(S(U)+\theta'\theta\eta)\cdot\eta\cdot \theta\eta\,
d\theta'd\theta,
\end{equation*}
or
\begin{equation*}
N(U,\eta)=\frac{1}{2} d^3F(S(U))\cdot\eta\cdot \eta
+\int_0^1\int_0^1 d^4F(S(U)+\theta''\theta'\theta\eta)\cdot
\theta'\theta\eta\cdot\eta\cdot \theta\eta\,
d\theta''d\theta'd\theta,
\end{equation*}
Then we use for~$s\in\R_+^*$,~$p,\,p_1,\,p_2,\,\in[1,\infty]$ such
that $\frac{1}{p}+\frac{s}{3}\geq\frac{1}{p_1}+\frac{1}{p_2}\geq
\frac{1}{p}$,

\begin{equation*}
\ds\|uv\|_{B^s_{p,q}}\leq C \|u\|_{B^{s}_{p_1,q}}
\|v\|_{B^s_{p_2,q}},
\end{equation*}
and for~$s>1$ , we use \cite[Proposition 2.1]{EscobedoVega}
\begin{equation*}
\ds\| d^{k}  F(\psi)\|_{B^s_{p_2,q}}\leq
C\left(s,F,\|\psi\|_{L^\infty}\right)\|\psi\|_{B^s_{p_2,q}},
\end{equation*}
for $k=3$ or $k=4$ and $d^{4}  F(z)=O(|z|)$, otherwise we decompose
$d^{4}  F(z)=A + O(|z|)$ where $A$ is a constant operator.

Eventually  using Lemma~\ref{lemma:ExpDecay} and
\begin{equation*}
\left\|\langle Q \rangle^\sigma \left|\eta\right|^l
\right\|_{B_{p_1,q}^s} \leq C\left\|\eta\right\|_{L^\infty}^{l-1}
\left\|\langle Q\rangle^\sigma\eta\right\|_{B_{p_1,q}^s},
\end{equation*}
for $l\in \N$, we conclude the proof.
\end{proof}

\begin{lemma}
                                                                \label{Lem:EstimateONF}
If Assumptions \ref{assumption:1}--\ref{assumption:NonLinearity} hold.
Let ~$\sigma\in\R$,~$s>1$, $p,p_1,p_2,q\in[1,\infty]$  and
$\sigma_1,\sigma_2\in\R$ such that
\begin{equation*}
\frac{1}{p}+\frac{s}{3}\geq\frac{1}{p_1}+\frac{1}{p_2}\geq\frac{1}{p}.
\end{equation*}
Then there exist~$\e>0$ and~$C>0$ such that for all~$U\in
B_\C(0,\e)$ and~$\eta \in B_{p,q}^s(\R^3,\R^8)\cap
L^\infty(\R^3,\R^8)$ with $\langle Q\rangle^{\sigma_1} \eta \in B_{p_1,q}^s(\R^3,\R^8)$
and $\langle Q\rangle^{\sigma_2} \eta \in B_{p_2,q}^s(\R^3,\R^8)$, we have
\begin{multline*}
\left\|<Q>^\sigma\left(\nabla F(S(U)+\eta) -\nabla
F(S(U)-\nabla F(\eta)\right)\right\|_{B^s_{p,q}}\\\leq
C(s,F,|U|+\left\|\eta\right\|_{L^\infty})\left(\left|U\right|+\left\|
<Q>^{\sigma_1}\eta\right\|_{B^s_{p_1,q}}\right)\left|U\right|\left\|
<Q>^{\sigma_2}\eta\right\|_{B^s_{p_2,q}}.
\end{multline*}
\end{lemma}
\begin{proof}
The proof is similar to the one of Lemma \ref{Lem:EstimateOnN}.
\end{proof}
\begin{lemma}
                                                                \label{lemma:EstimateOnN1}
If Assumptions \ref{assumption:1}--\ref{assumption:NonLinearity} hold.
Let ~$\sigma\in\R$,~$s>1$ and~$p,q\in[1,\infty]$ such that~$sp\geq
3$. Then there exist~$\e>0$ and~$C>0$ such that for all $U,\,U'\in
B_{\C^2}(0,\e)$ and~$\eta,\,\eta' \in B_{p,q}^s(\R^3,\R^8)$, we have
\begin{multline*}
\left\|\langle
Q\rangle^{\sigma}\left\{N(U,\eta)-N(U',\eta')\right\}\right\|_{B_{p,q}^s}
\leq
C\left(s,F,\left|U\right|+\left|U'\right|+\left\|\eta\right\|_{B_{p,q}^s}
+\left\|\eta'\right\|_{B_{p,q}^s}\right)\times\\
\times\Bigg\{ \left(\|\langle
Q\rangle^{\sigma_1}\eta\|_{B_{p,q}^s}+\|\langle
Q\rangle^{\sigma_1}\eta'\|_{B_{p,q}^s}\right)^2
\left(\left|U-U'\right|+\left\|\langle
Q\rangle^{\sigma_2}\left(\eta-\eta'\right)\right\|_{B_{p,q}^s}\right)\\
+\left(\left|U\right| +\left|U'\right|+\left\|\langle
Q\rangle^{\sigma_1'}\eta\right\|_{B_{p,q}^s} +\left\|\langle
Q\rangle^{\sigma_1'}\eta'\right\|_{B_{p,q}^s}\right)\times
\\\times\left(\|\langle Q\rangle^{\sigma_2'}\eta\|_{B_{p,q}^s}+\|\langle
Q\rangle^{\sigma_2'}\eta'\|_{B_{p,q}^s}\right)\|\langle
Q\rangle^{\sigma_3'}\left(\eta-\eta'\right)\|_{B_{p,q}^s}\Bigg\},
\end{multline*}
with~$2\sigma_1+\sigma_2=\sigma_1'+\sigma_2'+\sigma_3'=\sigma$ if
$<Q>^w\eta,\,<Q>^w\eta' \in B_{p,q}^s(\R^3,\R^8)$ for
$w\in\left\{\sigma_1,\sigma_2,\sigma_1',\sigma_2',\sigma_3'\right\}$.
\end{lemma}
\begin{proof}
Using the identity
\begin{equation*}
N(u,\eta)=\int_0^1\int_0^1
d^3F(S(u)+\theta'\theta\eta)\cdot\eta\cdot \theta\eta\,
d\theta'd\theta.
\end{equation*}
we can restrict the study to~$d^3 F(\phi)-d^3 F(\phi')$. If $F
=O(|z|^5)$, we have
\begin{equation*}
\|\langle Q\rangle^\sigma \left(d^3 F(\phi)-d^3
F(\phi')\right)\|_{B_{p,q}^s}\leq\int_0^1 \| d^{4}
F(\phi+t(\phi-\phi'))\|_{B_{p,q}^s}\|\langle
Q\rangle^\sigma(\phi-\phi')\|_{B_{p,q}^s}\,dt.
\end{equation*}
Then since~$s>1$ and $sp\geq 3$, we use
\begin{equation*}
\ds\|d^{4}  F(\psi)\|_{B^s_{p,q}}\leq C(s,F, \|\psi\|_{B^s_{p,q}}).
\end{equation*}
Using Lemma~\ref{lemma:ExpDecay}, we conclude the proof when $F
=O(|z|^5)$.

Otherwise, if $F$ is an homogeneous polynomial of order $4$, the
proof is easily adaptable since~$d^4 F$ is a constant tensor.

The case $F=O(|z|^4)$ follows by summing the two previous one since
as a function of $u\in\R^8$, $F(u)=Au^{\otimes4}+O(|u|^5) $.
\end{proof}

\begin{lemma}
                                                                \label{Lem:EstimateOnA}
If Assumptions \ref{assumption:1}--\ref{assumption:NonLinearity} hold.
Let ~$\sigma\in\R$,~$s\in\R$ and~$p,q\in[1,\infty]$. Then there
exist ~$\e>0$, ~$M>0$ and ~$C>0$ such that for all~$U,\, U'\in
B_{\C^2}(0,\e)$ and~$\eta,\, \eta' \in B_{L^2(\R^3,\R^8)}(0,M)$ with
$\langle Q\rangle^\sigma\left\{\eta-\eta'\right\}\in
B_{p,q}^s(\R^3,\R^8)$, one has
\begin{equation}
\left|A(U,\eta)-A(U',\eta')\right|\leq
C\left\{\left|U-U'\right|+\left\|\langle
Q\rangle^\sigma\left\{\eta-\eta'\right\}\right\|_{B^s_{p,q}}\right\}.
\end{equation}
\end{lemma}
\begin{proof}
We recall that
\begin{equation*}
A(U,\eta)=[\langle JdS(U),dS(U)\rangle - \langle
J\eta,d^2S(U)\rangle]^{-1}.
\end{equation*}
We have
\begin{multline*}
A(U,\eta)-A(u',\eta')=-[\langle JdS(U),dS(U)\rangle - \langle
J\eta,d^2S(U)\rangle]^{-1}\times\\
\times\left\{\langle JdS(U),dS(U)\rangle - \langle
J\eta,d^2S(U)\rangle-\langle JdS(U'),dS(U')\rangle + \langle
J\eta',d^2S(U')\rangle\right\}\times\\
\times[\langle JdS(U'),dS(U')\rangle - \langle
J\eta',d^2S(U')\rangle]^{-1}.
\end{multline*}
The lemma then follows from Lemma~\ref{lemma:ExpDecay}.
\end{proof}

%
%
\subsubsection{Global wellposedness for~$z$}

Let ~$U\in {\mathcal U}(\e,\delta)$ and~$z_0\in {\mathcal
H}_c(U(0))\cap H^{s}$. Let us write $U_\infty=\ds\lim_{t\to +\infty}
U(t)$, we define~${\mathcal T}_{U,z_0}(z)$ by
\begin{multline*}
{\mathcal T}_{U,z_0}(z)(t)= e^{-\rmi tH+\rmi\int_0^tE(U(r))\;dr}z_0
+\int_0^t e^{-\rmi(t-v)H+\rmi\int_v^tE(U(r))\;dr}{\mathbf P}_c(U(v))J\nabla F(\eta(v))\,dv\\
+\int_0^te^{-\rmi(t-v)H+\rmi\int_v^tE(U(r))\;dr}{\mathbf
P}_c(U(v))J\{\nabla F(S(U(v))+\eta(v))
-\nabla F(S(U(v))-\nabla F(\eta(v))\}\,dv\\
+\int_0^te^{-\rmi(t-v)H+\rmi\int_v^tE(U(r))\;dr}{\mathbf
P}_c(U(v))dS(U(v))A(U(v),\eta(v))
\langle N(U(v),\eta(v)) ,dS(U(v))\rangle\,dv\\
-\int_0^te^{-\rmi(t-v)H+\rmi\int_v^tE(U(r))\;dr}(d {\mathbf
P}_c(U(v))\dot{U}(v)\eta(v)\,dv.
\end{multline*}
with
\begin{equation*}
\eta(t)=z(t) +g\left(U(t),z(t)\right)
\end{equation*}
First, we have a local wellposedness result for $z$ with the
\begin{lemma}
                                                                \label{Lem:LocalWellPosednessForZ}
If Assumptions \ref{assumption:1}--\ref{assumption:Resonance} hold.
Then there exist~$\e_0>0$ and~$\delta_0>0$ such that for any
$\e\in(0,\e_0)$, $\delta\in(0,\delta_0)$, ~$U\in {\mathcal
U}(\delta,\e)$ and $z_0\in B_{H^{s}(0,\delta)}\cap {\mathcal
H}_c(U(0))$ there are $T^\pm(z_0,U)>0$ and a solution
\begin{equation*}
\ds z\in \cap_{k=0}^2 {\mathcal
C}^k((-T^-(z_0,U);+T^+(z_0,U)),H^{s-k}(0,\delta))
\end{equation*}
of the equation
\begin{equation}
                                                                \label{Eq:ForZ}
\begin{cases} \p_t z=JH(U)z+{\mathbf P}_c(U)JN(u,\eta)-
( d{\mathbf P}_c(U))\dot{U}\eta,\\
z(0)=z_0,
\end{cases}
\end{equation}
where~$\eta(t)=z(t) +g\left(U(t),z(t)\right)$.

Moreover, $z$ is unique in $L^\infty((-T',T),H^{s})$ for any $T\in
(0,T^+(z_0,U))$ and $T'\in (0,T^-(z_0,U))$ and we have if
$T^+(z_0,U)<+\infty$ then $$\ds\lim_{t\to
T^+(z_0,U)}\left\|z(t)\right\|_{H^{s}} \geq \delta$$ and if 
$T^-(z_0,U)=+\infty$ then $$\ds\lim_{t\to -T^-(z_0,U)}
\left\|z(t)\right\|_{H^{s}}\geq~\delta.$$
\end{lemma}
\begin{proof}
It is a consequence of the fix point theorem applied to~${\mathcal
T}_{U,z_0}$:

Using Lemmas \ref{Lem:EstimateOnN}, \ref{lemma:EstimateOnN1} and
\ref{Lem:EstimateOnA} with the Estimate
\eqref{Estimate:NormLinearizedFlow1}--\eqref{Estimate:NormLinearizedFlow3}
and the properties of $g$ given by Lemma
\ref{Lem:DecompositionCentralManifofd2}, we obtain that ${\mathcal
T}_{U,z_0}$ leaves a small ball in $H^{s}$ invariant and is a
contraction inside this ball.

Hence there exists a unique solution defined on the
interval~$[-T,T]$. Classical arguments permit to extend the solution
over a maximal interval $(-T^-(z_0,U),T^+(z_0,U))$ such that
if~$T^+(z_0,U)<\infty$ then necessarily the solution should leave a
small ball in~$H^{s}$  at time $T^+(z_0,U)$.
\end{proof}

We have now a global wellposedness result as stated in the
\begin{lemma}
                                                                \label{Lem:StabilizationForT}
If Assumptions \ref{assumption:1}--\ref{assumption:Resonance} hold.
There exist~$\e_0>0$, $\delta_0>0$ and~$C>0$ such that for any
$\e\in(0,\e_0)$, $\delta\in(0,\delta_0)$, $U\in {\mathcal
U}(\e,\delta)$ and $z_0\in B_{H^{s}}(0,\delta)\cap
\mathcal{H}_c(U(0))$ we obtain for the Cauchy problem
\eqref{Eq:ForZ}, $T^+(U,z_0)=+\infty$, $T^-(U,z_0)=+\infty$,
$z\in{\mathcal Z}(U,\delta)$ and
\begin{equation*}
\max\left[\|z\|_{L^\infty(\R,H^{s})},
\|z\|_{L^2(\R,H^{s}_{-\sigma})},
\|z\|_{L^2(\R,B^\beta_{\infty,2})}\right]\leq
C\left\|z_0\right\|_{H^{s}}.
\end{equation*}
\end{lemma}
\begin{proof}
We have $(1-P_c(U))z \equiv 0$ because its time derivative is zero
and $(1-P_c(U(0)))z(0)=0$.

Let us introduce for any $0<T< T^+(U,z_0)$, the function
\begin{equation*}
    m(T)=\sup_{t\in(-T,T)}\left\{\left\|z\right\|_{L^\infty((-T,T),H^{s})},\;
    \left\|z\right\|_{L^2((-T,T),H^{s}_{-\sigma})},\;\|z\|_{L^2((-T,T),B^\beta_{\infty,2})}\right\}
\end{equation*}

First, we study the estimation of $L^2((-T,T),H^{s}_{-\sigma})$ . We
use Estimate \eqref{Eq:LowerBoundPc} and the estimates of the
Theorem B.1.
\begin{eqnarray*}
\lefteqn{\left\|z\right\|_{L^2((-T,T),H^{s}_{-\sigma})}}\\
&\leq&C_0\left\|z_0\right\|_{H^{s}} +C\left\|{\mathbf P}_c\int_0^t
e^{-\rmi(t-v)H+\rmi\int_v^tE(U(r))\;dr}
{\mathbf P}_c(U(v))J\nabla F(\eta(v))\,dv\right\|_{L^2((-T,T),H^{s}_{-\sigma})}\\
&&+C\left\|\nabla F(S(U)+\eta)
-\nabla F(S(U)-\nabla F(\eta)\right\|_{L^2((-T,T),H^{s}_{\sigma})}\\
&&+C\left\|dS(U)A(U,\eta)
\langle N(U,\eta) ,dS(U)\rangle\right\|_{L^2((-T,T),H^{s}_{\sigma})}\\
&&+C\left\|(d {\mathbf
P}_c(U)\dot{U}\eta\right\|_{L^2((-T,T),H^{s}_{\sigma})}.
\end{eqnarray*}
We now study the estimation of the third term of the right hand side
\begin{eqnarray*}
    \lefteqn{\left\|\int_0^t e^{-\rmi(t-v)H+\rmi\int_v^tE(U(r))\;dr}{\mathbf P}_c{\mathbf P}_c(U(v))J\nabla
    F(\eta(v))\,dv\right\|_{L^2_t((-T,T),H^{s}_{-\sigma})}}\\
    &\leq& \int_{-T}^T \left\|e^{-\rmi(t-v)H+\rmi\int_v^tE(U(r))\;dr}{\mathbf P}_c{\mathbf P}_c(U(v))J\nabla
    F(\eta(v))\right\|_{L^2_t((-T,T),H^{s}_{-\sigma})}\,dv\\
    &\leq& C(U)\left\|\nabla F(\eta)\right\|_{L^1((-T,T),H^{s})}\qquad\qquad\qquad\qquad\qquad\qquad\qquad\qquad\qquad\\
    &\leq&
    C(U)\left\|\eta\right\|_{L^2((-T,T),L^\infty)}^2\left\|\eta\right\|_{L^\infty((-T,T),H^{s})},
\end{eqnarray*}
where we used Theorem B.1 Estimate $(ii)$.
Hence for the $L^2H^{s}_{-\sigma}$ estimate, we obtain
\begin{eqnarray*}
\left\|z\right\|_{L^2((-T,T),H^{s}_{-\sigma})}
&\leq&C_0\left\|z_0\right\|_{H^{s}}
+C\left\|\eta\right\|_{L^2((-T,T),L^\infty)}^2\left\|\eta\right\|_{L^\infty((-T,T),H^{s})}\\
&&+C\left(\left\|U\right\|_{\infty}+\left\|\eta\right\|_{L^\infty((-T,T),H^{s}_{-\sigma})}\right)\left\|U\right\|_{\infty}\left\|
\eta\right\|_{L^2((-T,T),H^{s}_{-\sigma})}\\
&&+C\left\|\eta\right\|_{L^2((-T,T),L^\infty)}^2+C\left\|\dot{U}\right\|_{L^2}\left\|\eta\right\|_{L^\infty((-T,T),H^{s})},
\end{eqnarray*}
using Lemma \ref{Lem:DecompositionCentralManifofd2}, we obtain
\begin{equation*}
\left\|z\right\|_{L^2((-T,T),H^{s}_{-\sigma})}\leq
C_0\left\|z_0\right\|_{H^{s}} +Cm(T)^3 +Cm(T)^2+ C\e m(T)+C\delta^2
m(T),
\end{equation*}
where $C$ depends of $\left\|U\right\|_{\infty}$ and
$\left\|\eta\right\|_{L^\infty((-T,T),H^{s})}$.
\medskip

Then, we estimate the $H^{s}$ norm. Using Estimate
\eqref{Eq:LowerBoundPc}, we have
\begin{eqnarray*}
\lefteqn{\left\|z(t)\right\|_{H^{s}} \leq \left\|z_0\right\|_{H^{s}}
+\int_{-T}^T \left\|\nabla F(\eta(v))\right\|_{H^{s}}\,dv} \\
&&+\Big\|\int_0^te^{-\rmi(t-v)H+\rmi\int_v^tE(U(r))\;dr}{\mathbf
P}_c(U(v))\times\\
&&\quad\times J\{\nabla F(S(U(v))+\eta(v))
-\nabla F(S(U(v))-\nabla F(\eta(v))\}\,dv\big\|_{H^{s}} \\
&&+\int_{-T}^T\left\|dS(U(v))A(U(v),\eta(v))
\langle N(U(v),\eta(v)) ,dS(U(v))\rangle\right\|_{H^{s}}\,dv \\
&&+\int_{-T}^T\left\|(d {\mathbf
P}_c(U(v))\dot{U}(v)\eta(v)\right\|_{H^{s}}\,dv.
\end{eqnarray*}
to estimate the third term of the right hand side, we use the
$H$-smoothness estimates, more precisely Theorem B.1 Estimate $(ii)$
and then we use Lemma B.14:
\begin{eqnarray*}
\lefteqn{\left\|\int_0^te^{-\rmi(t-v)H+\rmi\int_v^tE(U(r))\;dr}{\mathbf
P}_c(U(v))J\{\nabla F(S(U(v))+\eta(v)) -\nabla F(S(U(v)))-\nabla
F(\eta(v))\}\,dv\right\|_{H^{s}}}\qquad\qquad\qquad\qquad\qquad\qquad\\
\lefteqn{\leq \left\|\int_0^te^{\rmi
vH-\rmi\int_0^vE(U(r))\;dr}{\mathbf P}_c(U(v))J\{\nabla
F(S(U(v))+\eta(v)) -\nabla F(S(U(v)))-\nabla
F(\eta(v))\}\,dv\right\|_{H^{s}}}\qquad\qquad\qquad\qquad\qquad\qquad\\
&&\leq C\left\|\{\nabla F(S(U)+\eta) -\nabla F(S(U)-\nabla
F(\eta)\right\|_{L^2((-T,T),H_\sigma^{s})}\\
&&\leq
C\left(\left\|U\right\|_{\infty}+\left\|\eta(v)\right\|_{L^\infty((-T,T),H^{s})}\right)\left\|U\right\|_{\infty}\left\|
\eta\right\|_{L^2((-T,T),H^{s}_{-\sigma})}
\end{eqnarray*}
Hence for the $L^\infty H^{s}$ estimate, we obtain
\begin{multline*}
\left\|z(t)\right\|_{H^{s}} \leq \left\|z_0\right\|_{H^{s}}
+C\left\|\eta\right\|_{L^\infty((-T,T), H^{s})}\left\|\eta\right\|_{L^2((-T,T), L^{\infty})}^2 \\
C\left(\left\|U\right\|_{L^\infty((-T,T))}+\left\|\eta(v)\right\|_{L^\infty((-T,T),H^{s})}\right)\left\|U\right\|_{L^\infty((-T,T))}\left\|
\eta\right\|_{L^2((-T,T),H^{s}_{-\sigma})} \\
+C\left\|\eta\right\|_{L^2((-T,T), L^{\infty})}^2
+\left\|\dot{U}\right\|_{L^1((-T,T))}\left\|\eta\right\|_{L^\infty((-T,T),H^{s})},
\end{multline*}
using Lemma \ref{Lem:DecompositionCentralManifofd2}, we obtain
\begin{equation*}
\left\|z(t)\right\|_{H^{s}} \leq \left\|z_0\right\|_{H^{s}} +Cm(T)^3
+ Cm(T)^2+C\e m(T)+ C\delta^2 m(T),
\end{equation*}
where $C$ depends of $\left\|U\right\|_{\infty}$ and
$\left\|\eta\right\|_{L^\infty((-T,T),H^{s})}$.

For the $L^2B^{\beta}_{\infty,2}$ estimate, by
Proposition~\ref{Prop:ContSpectrumLinOp} and Theorem
\ref{Thm:Strichartz}, we have for any $\e>0$, any $p_\e>3/\e$ and
$\theta_{\e}=\frac{4}{p_\e-2}$
\begin{eqnarray*}
\|z\|_{L^2((-T,T),B^\beta_{\infty,2})}&\leq& \|z\|_{L^2((-T,T),B^{\beta+\e}_{p_\e,2})}\\
&\leq&C_0\|z_0\|_{H^{\beta+1+\e+\theta_\e/2}}
+C\left\|d^2F(S(U))\cdot\eta\right\|_{L^2((-T,T),B^{\beta+2+\e+\theta_\e}_{p_\e',2})}\\
&&+C\left\|N(U,\eta)\right\|_{L^1((-T,T),H^{\beta+1+\e+\theta_\e/2})}\\
&&+C\left\|dS(U)A(U,\eta) \langle N(U,\eta),dS(U)\rangle\right\|_{L^1((-T,T),H^{\beta+1+\e+\theta_\e/2})}\\
&&+C\left\|(d {\mathbf
P}_c(U))\dot{U}\eta\right\|_{L^1((-T,T),H^{\beta+1+\e+\theta_\e/2})}\,dv.
\end{eqnarray*}
With Lemma~\ref{Lem:EstimateOnN} and \ref{Lem:EstimateONF}, we infer
\begin{eqnarray*}
\lefteqn{\|z\|_{L^2(\R,B^\beta_{\infty,2})}\leq C_0\|z_0\|_{H^{\beta+1+\e+\theta_\e/2}}+C|U|_{\infty}\left\|\eta\right\|_{L^2((-T,T),H_{-\sigma}^{\beta+2+\e+\theta_\e})}}\\
&&+C|U|_\infty\left\|z\right\|_{L^2((-T,T),L^{\infty})}\left\|z\right\|_{L^2((-T,T),H^{\beta+1+\e+\theta_\e/2})}
\\
&&+C(|U|_\infty+
\left\|\eta\right\|_{L^\infty((-T,T),H^{\beta+1+\e+\theta_\e/2})})
\left\|\eta\right\|_{L^2((-T,T),L^\infty)}^2\left\|z\right\|_{L^\infty((-T,T),H^{\beta+1+\e+\theta_\e/2})}\\
&&+C(|U|_\infty+
\left\|\eta\right\|_{L^\infty((-T,T),H^{\beta+1+\e+\theta_\e/2})})
\left\|\eta\right\|_{L^2((-T,T),H_{-\sigma}^{\beta+1+\e+\theta_\e/2})}\left\|\eta\right\|_{L^\infty((-T,T),H^{\beta+1+\e+\theta_\e/2})}\\
&&+C\|\dot{U}\|_{L^1}\left\|\eta\right\|_{L^\infty((-T,T),H^{\beta+1+\e+\theta_\e/2})},
\end{eqnarray*}
we infer since for small $\e>0$, $s\geq\beta+2+\e+\theta_\e$ and
using Lemma \ref{Lem:DecompositionCentralManifofd2},
\begin{eqnarray*}
\|z\|_{L^2((-T,T),B^{\beta}_{\infty,2})}
&\leq&C_0\|z_0\|_{H^{\beta+1+\e+\theta_\e/2}} +Cm(T)^3 +Cm(T)^2+C\e
m(T)+C\delta^2m(T).
\end{eqnarray*}
Hence we obtain
\begin{eqnarray*}
m(T) &\leq&C_0\|z_0\|_{H^{\beta+1+\e+\theta_\e/2}}+C\e
m(T)+C\delta^2 m(T)+Cm(T)^3 +Cm(T)^2,
\end{eqnarray*}
where $C_0$ do not depend of $m$ and $C$ is a nondecreasing function
of $\left\|z\right\|_{L^\infty((-T,T),H^s)}$ and
$\left\|U\right\|_\infty$ and hence it can be bounded by a
nondecreasing function of $m$.

If $\|z_0\|_{H^{s}}$ is small then $m(0)$ is small and $m(T)$ stay
small. Therefore we have that $z\in{\mathcal Z}(U,\delta)$
if~$\|z_0\|_{H^{s}}$is small enough for any $\delta$ and $\e$ are
small enough and
\begin{equation*}
\max\left[\|z\|_{L^\infty(\R,H^{s})},
\|z\|_{L^2(\R,H^{s}_{-\sigma})},
\|z\|_{L^2(\R,B^\beta_{\infty,2})}\right]\leq
f(\left\|z_0\right\|_{H^{s}})
\end{equation*}
where $f$ is such that there exists $C>0$ with
\begin{equation*}
f(\left\|z_0\right\|_{H^{s}})\leq C\left\|z_0\right\|_{H^{s}}.
\end{equation*}
\end{proof}

The solution $z$ just found is a function of $z_0$ and $U$, writing
it $z[z_0,U]$, we have the following important property given by the
\begin{lemma}
                                                                \label{Lem:LipshitzPropertyOfT}
If Assumptions \ref{assumption:1}--\ref{assumption:Resonance} hold.
Then for any $T>0$, there exist~$\e_0>0$, ~$\delta_0>0$, $C>0$ and
$\kappa\in(0,1)$ such that for any $\e\in(0,\e_0)$,
$\delta\in(0,\delta_0)$, $U,\,U'\in {\mathcal U}(\e,\delta)$,
$z_0\in {\mathcal H}_c(U(0))$, $z_0'\in {\mathcal H}_c(U'(0))$,
$z\in {\mathcal Z}(U,\delta)$ and $z'\in {\mathcal Z}(U',\delta)$,
one has
\begin{multline*}
\left\|z[z_0',U']
-z[z_0,U]\right\|_{L^\infty((-T;T),H^{s})\cap L^2((-T;T),L^{\infty})\cap L^2((-T;T),H^{s}_{-\sigma})}\\
\leq C \left\|z_0-z_0'\right\|_{H^{s}} +\kappa\left\{
\left\|U-U'\right\|_{L^\infty((-T;T))}+\left\|\dot{U}-\dot{U}'\right\|_{L^\infty((-T;T))}\right\}.
\end{multline*}
\end{lemma}
\begin{proof}
We use the technics of the previous lemma.
\end{proof}


\subsubsection{Global wellposedness for~$U$ and its stabilization}

Here we want to solve the equation for~$U$. We notice that~$z$ and
$\alpha$ have been built in the previous section and are functions
of~$U$ and~$z_0\in {\mathcal H}_c(U(0))$. Let us introduce for any
~$U_0\in B_\C(0,\e)$ the function on ${\mathcal U}(\e,\delta)$:
\begin{equation*}
f_{U_0}(U)(t)=U_0-\int_0^tA(U(v),\eta(v))\langle N(U(v),\eta(v))
,dS(U(v))\rangle\,dv,
\end{equation*}
where~$\eta=z(t) +g\left[U(t),z(t)\right]$. We have the
\begin{lemma}
If Assumptions \ref{assumption:1}--\ref{assumption:Resonance} hold.
There exist~$\e_0>0$ and~$\delta_0>0$ such that for any
$\e\in(0,\e_0)$, $\delta\in(0,\delta_0)$, the function $f_{U_0}$
maps~${\mathcal U}(\e,\delta)$ into itself if~$U_0$ and $z_0\in
H^{s}\cap {\mathcal H}_c(U_0)$ are small enough.
\end{lemma}
\begin{proof}
By means of Lemma~\ref{Lem:EstimateOnN}, we obtain
\begin{equation*}
\|\p_tf_{U_0}(U)\|_{L^1(\R)\cap L^\infty(\R) }\leq
C\left\|N(U(v),\eta(v))\right\|_{L^1(\R,H^{s}_{-\sigma})\cap
L^\infty (\R,H^{s})} \leq
 \delta^2.
\end{equation*}
and
\begin{equation*}
\|f_{U_0}(U)\|_{L^\infty(\R)}\leq |U_0|+
C\left\|N(U(v),\eta(v))\right\|_{L^1(\R,H^{s})} \leq |U_0|+
 \delta^2,
\end{equation*}
hence for sufficiently small $U_0$ and $\delta$, we obtain the
lemma.
\end{proof}

The function~$f_{U_0}$ has also a local Lipshitz property as stated by the
\begin{lemma}                                                   \label{Lem:LipshitzPropertyOff}
If Assumptions \ref{assumption:1}--\ref{assumption:Resonance} hold.
For any $T>0$, there exist~$\e_0>0$,~$\delta_0>0$ and
$\kappa\in(0,1)$ such that for any $\e\in(0,\e_0)$,
$\delta\in(0,\delta_0)$, ~$U,\,U'\in {\mathcal U}(\e,\delta)$, for
any~$z_0\in {\mathcal H}_c(U(0))\cap H^{s}$, for any $z_0'\in
{\mathcal H}_c(U'(0))\cap H^{s}$ small enough, for $U_0,U_0'$ small
enough, such that
\begin{multline*}
\left|f_{U_0}(U)-f_{U_0'}(U')\right|_{L^\infty((-T;T))}
+\left|\p_tf_{U_0}(U)-\p_tf_{U_0'}(U')\right|_{L^1((-T;T))}\\
\leq \left|U_0-U_0'\right|+
\kappa\left(\left\|U-U'\right\|_{L^\infty((-T;T))}
+\left\|\dot{U}-\dot{U}'\right\|_{L^1((-T;T))}+\|z_0-z_0'\|_{H^{s}}\right).
\end{multline*}
\end{lemma}
\begin{proof}
This a straightforward consequence of
Lemma~\ref{lemma:EstimateOnN1},~\ref{Lem:EstimateOnA}
and~\ref{Lem:LipshitzPropertyOfT}.
\end{proof}

We now obtain the
\begin{lemma}
If Assumptions \ref{assumption:1}--\ref{assumption:Resonance} hold.
There exists~$\e>0$ and $\delta>0$ such that for any~$U_0\in \C$
small and~$z_0\in {\mathcal H}_c(U_0)\cap H^{s}_\sigma$ small, the
equation
\begin{equation}                                                \label{Eq:forU}
\left\{\begin{array}{lll}
\dot{U}&=&-A(U,\eta)\langle N(U,\eta) ,dS(U)\rangle,\\
U(0)&=&U_0,
\end{array}\right.
\end{equation}
where~$\eta(t)=z(t) +g\left[U(t),z(t)\right]$, has a unique solution
in~${\mathcal U}(\delta,\e)$. Moreover, there exists $C>0$ such that
\begin{equation*}
    \left|U_{\pm\infty}-U_0\right|\leq C
\left\|z_0\right\|_{H^{s}}^2.
\end{equation*}

\end{lemma}
\begin{proof}
This is also a the fixed point
result for~$f_{U_0}$. Let us fix $T>0$ and consider, for any $V\in {\cal U}(\delta,\e)$
with sufficiently small $\delta>0$ and $\e>0$, the sequence:
\begin{equation*}
\begin{cases}
    V_{n+1}=f_{U_0}(V_n),\; \forall n\in\N\\
    V_0=V;
\end{cases}
\end{equation*}
for any $n\in \N$, $V_n\in  {\cal U}(\delta,\e)$. With Lemma
\ref{Lem:LipshitzPropertyOff}, the fixed point theorem give us the
convergence for the norms of $L^\infty((-T,T))$ and
$\dot{W}^{1,1}((-T,T))$ of $\left(V_n\right)_{n\in\N}$.

Then we notice that for any $T'\in\R$, we have
\begin{equation*}
        V_{n+1}(t)=f_{f_{U_0}(V_n)(T')}(V_n)(t-T').
\end{equation*}
Since for $T'\in(-T;T)$, $\left(f_{U_0}(V_n)(T')\right)$ is a Cauchy
sequence, the Lemma \ref{Lem:LipshitzPropertyOff} give us the
convergence of $\left(V_n\right)$ for the norms of
$L^\infty((T'-T;T'+T))$ and $\dot{W}^{1,1}((T'-T;T'+T))$.

Iterating this process, we obtain that the sequence converges uniformly 
locally in time and we prove the lemma since the other
statements are classical. We just notice that the last statement
follow from the fact that there exists $C>0$ such that
\begin{equation*}
    \int_{\R^\pm} \left|\dot{U}(v)\right|\;dv\leq \int_{\R^\pm}
\left|A(U(v),\eta(v))\langle N(U(v),\eta(v)) ,dS(U(v))\rangle
\right|\;dv \leq C \left\|z_0\right\|_{H^{s}}^2.
\end{equation*}
\end{proof}





\subsubsection{The asymptotic profile of~$z$}
In this section, our aim is to precise the asymptotic profile of
$z$ when $z_0$ is localized. 
First we state the
\begin{proposition}                                             \label{Prop:WeightedEstimates}
There exists $\e>0$, such that for
all $U\in B_{\C^2}(0,\e)$ and $\alpha\in \R^+$ there exists $C>0$
such that
\begin{equation*}
    \left\|\left\langle Q\right\rangle^\alpha e^{J t
H(U)}\psi\right\|\leq C_\alpha \sum_{\beta=0}^\alpha\left\langle
t\right\rangle^\beta\left\|\left\langle
Q\right\rangle^{\alpha-\beta}\psi\right\|
\end{equation*}
for any $\psi\in L^2(\R^3,\C^8)$.
\end{proposition}
\begin{proof}
From Proposition \ref{Prop:ContSpectrumLinOp}, we obtain the result
for $\alpha=0$, then we just need the result the estimate
\begin{equation*}
    \left\|Q^\alpha e^{J t
H(U)}\psi\right\|^2\leq C_\alpha^2 \sum_{0\leq \beta \leq\alpha}\left|
t\right|^{2|\beta|}\left\| Q^{\alpha-\beta}\psi\right\|^2
\end{equation*} for any
$\psi\in L^2(\R^3,\C^8)$, $\alpha\in \N^3$ and some $C>0$
independent of $\psi$. The rest of the proposition will follow by
interpolation.

For $U=0$, this follows by an iterated proof from the identity
\begin{equation*}
    \frac{d}{dt}e^{ \rmi t H}Qe^{-\rmi t H}=e^{ \rmi t H}\alpha e^{-\rmi t H}
\end{equation*}
where $\alpha$ is the $3$-vector of Dirac Pauli matrices defined in
the introduction. For $U\neq 0$, we use the same proof with the
exponential decay of Proposition \ref{Prop:ManifoldPLS}.
\end{proof}

We can improve Lemma \ref{Lem:StabilizationForTloc}, if we use
\cite[Theorem 1.2]{NBStableDirectionSmallSolitonNR} and
\cite[Theorem 1.1]{NBStableDirectionSmallSolitonNR}, which we repeat
here :
\begin{theorem}[Theorem 1.1 of \cite{NBStableDirectionSmallSolitonNR}:
 Propagation for perturbed Dirac dynamics]
                                                                \label{Thm:Propagation}
Assume that Assumptions~\ref{assumption:1} and~\ref{assumption:2}
hold and let be~$\sigma>\frac{5}{2}$. Then one has
\begin{equation*}
\|e^{-\rmi t H}{\mathbf P}_c\left(H\right)\|_{
B(L^2_{\sigma},L^2_{-\sigma})}\leq C\left\langle
t\right\rangle^{-\frac{3}{2}}.
\end{equation*}
\end{theorem}
We also have
\begin{proposition}[Proposition 2.2 of 
\cite{NBStableDirectionSmallSolitonNR}: Propagation far from thresholds]
                                                                   \label{Prop:PropagationHighEnergy}
Suppose that Assumption~\ref{assumption:1} holds. Then for any
$\chi\in {\mathcal C}^\infty(\R^3,\C^4)$ bounded with support in
$\R\setminus(-m;m)$ and for any~$\sigma\geq 0$, there is $C>0$ such
that
\begin{equation*}
\|e^{-\rmi t H}\chi\left(H\right)\|_{
B(L^2_{\sigma},L^2_{-\sigma})}\leq C\left\langle
t\right\rangle^{-\sigma}.
\end{equation*}
\end{proposition}
Using Duhamel's formula like in Proposition \ref{Prop:ContSpectrumLinOp} 
and interpolating with estimate \eqref{Eq:PertutbedConservationLaw}, we obtain the
\begin{corrollary}																																			\label{Cor:LinPropagation}
Assume that Assumptions~\ref{assumption:1} and~\ref{assumption:2}
hold and let $\theta \geq 0$ and ~$\sigma>\frac{5}{2}\theta$. Then there exists $\e>0$, such that for
all $U\in B_{\C^2}(0,\e)$ one has
\begin{equation*}
\|e^{J t H(U)}{\mathbf P}_c\left(U\right)\|_{
B(L^2_{\sigma},L^2_{-\sigma})}\leq C\left\langle
t\right\rangle^{-\frac{3\theta}{2}}.
\end{equation*}
\end{corrollary}

\begin{theorem}[Theorem 1.1 of \cite{NBStableDirectionSmallSolitonNR}:
Dispersion for perturbed Dirac dynamics]
                                                               \label{Thm:Dispersion}
Assume that Assumptions~\ref{assumption:1} and~\ref{assumption:2}
hold. Then for~$p\in[1,2]$,~$\theta\in[0,1]$, $s-s'\geq (2
+\theta)(\frac{2}{p}-1)$ and $q\in[1,\infty]$ there exists a
constant~$C>0$ such that
\begin{equation*}
\|e^{-\rmi t H}{\mathbf P}_c(H)\|_{B^s_{p,q},B^{s'}_{p',q}} \leq C
\left(K(t)\right)^{\frac{2}{p}-1}
\end{equation*}
with~$\frac{1}{p}+\frac{1}{p'}=1$, and
\begin{equation*}
K(t)= \left\{
\begin{array}{ll}
\ds\left|  t\right| ^{-1+\theta/2}
& \mbox{if } |t|\in (0,1],\vspace{0.3cm}\\
\ds\left|  t\right| ^{-1-\theta/2} & \mbox{if } |t|\in [1,\infty).
\end{array}
\right.
\end{equation*}
\end{theorem}

Using, once more Duhamel's formula, the previous theorem and corollary \ref{Cor:LinPropagation}, we obtain the 
\begin{corrollary}																																			\label{Cor:LinDispersion}
Assume that Assumptions~\ref{assumption:1} and~\ref{assumption:2}
hold and let be ~$p\in[1,2]$,~$\theta\in[0,1]$, $s-s'\geq (2
+\theta)(\frac{2}{p}-1)$, $q\in[1,\infty]$ and ~$\sigma>\max\{\frac{3}{2}, (\frac{2}{p}-1)(1+\frac{\theta}{2})\}$. Then there exists $\e>0$, such that for
all $U\in B_{\C^2}(0,\e)$ one has
\begin{equation*}
\|e^{J t H(U)}{\mathbf P}_c(U)\|_{H^s_{\sigma},B^{s'}_{p',q}} \leq C
\left(K(t)\right)^{\frac{2}{p}-1}
\end{equation*}
with~$\frac{1}{p}+\frac{1}{p'}=1$, and
\begin{equation*}
K(t)= \left\{
\begin{array}{ll}
\ds\left|  t\right| ^{-1+\theta/2}
& \mbox{if } |t|\in (0,1],\vspace{0.3cm}\\
\ds\left|  t\right| ^{-1-\theta/2} & \mbox{if } |t|\in [1,\infty).
\end{array}
\right.
\end{equation*}
\end{corrollary}
\begin{proof}
We first prove it for $U=0$. We have to study the high and low energy part in a different manner. 
For the low energy part, we iterate twice Duhamel's formula with respect to $D_m$ in order to use 
Theorem \ref{Thm:Propagation} and  Theorem \ref{Thm:Dispersion} for the free case.

In the high energy part, we use also Duhamel's formula. But, we use Poposition 
Theorem \ref{Thm:Dispersion} for the free case and Proposition \ref{Prop:PropagationHighEnergy}. 

Then for $U\neq 0$, we work like for Estimate \eqref{Eq:PertutbedConservationLaw}. 
\end{proof}

We obtain the
\begin{lemma}
                                                                \label{Lem:StabilizationForTloc}
If Assumptions \ref{assumption:1}--\ref{assumption:Resonance} hold.
There exist~$\e_0>0$, $\delta_0>0$ and~$C>0$ such that for any
$\e\in(0,\e_0)$, $\delta\in(0,\delta_0)$, $U_0\in B_{\C^2}(0,\e)$
and $z_0\in B_{H^{s}_\sigma}(0,\delta)\cap \mathcal{H}_c(U_0)$ we
obtain for the Cauchy problem \eqref{Eq:ForZ} (with $U$ the solution
of \eqref{Eq:forU}) a global solution $z$ such that
\begin{equation*}
\max\left[\sup_{t\in \R}(\|z(t)\|_{H^{s}}),\; \sup_{t\in \R}(\langle
t\rangle^{3/2}\|z(t)\|_{H^{s}_{-\sigma}}),\; \sup_{t\in \R}\langle
t\rangle^{3/2}\|z(t)\|_{B^\beta_{\infty,2}}),\; \sup_{t\in
\R}(\langle t\rangle^{-3/2}\|z(t)\|_{H^{s}_{3/2}})\right]\leq
C\left\|z_0\right\|_{H^{s}_\sigma}.
\end{equation*}
\end{lemma}
\begin{proof}
The proof is similar to the one of Lemma \ref{Lem:StabilizationForT}
with some adaptations involving the norm 
$H^{s}_{\sigma}$, we also refer to the proof of \cite[Lemma
5.5]{NBStableDirectionSmallSolitonNR}.

Let
\begin{equation*}
t\mapsto \xi_\pm(t)=e^{J \int_0^t
\left(E(U(v))-E(U_{\pm\infty})\right)\;dv}z(t)
\end{equation*}
and
\begin{equation*}
t\mapsto V_\pm(t)=e^{-\rmi \int_0^t
\left(E(U(v))-E(U_{\pm\infty})\right)\;dv}U(t).
\end{equation*}
We use exactly the same method as the one of Lemma
\ref{Lem:StabilizationForT}, applied to
\begin{multline*}
\xi_\pm(t)=e^{Jt
H(V_{\pm\infty})}z_0 +\int_0^te^{J(t-s)
H(V_{\pm\infty})}{\mathbf P}_c(V_\pm(v))J\left(d^2 F(S(V_\pm(v)))
-d^2 F(S(V_{\pm\infty}))\right)\xi_\pm(v)\,dv\\
+\int_0^te^{J(t-s) H(V_{\pm\infty})}{\mathbf P}_c(V_\pm(v))JN(V_\pm(v),
\tilde{\eta}_\pm(v))\,dv\\
+\int_0^te^{J(t-s) H(V_{\pm\infty})}{\mathbf
P}_c(V_\pm(v))dS(V(v))A(V_\pm(v),\tilde{\eta}_\pm(v))
\langle N(V_\pm(v),\tilde{\eta}_\pm(v)) ,dS(V_\pm(v))\rangle\,dv\\
-\int_0^te^{J(t-s) H(V_{\pm\infty})}(d {\mathbf
P}_c(V_\pm(v)))A(V_\pm(v),\tilde{\eta}_\pm(v))\langle
N(V_\pm(v),\tilde{\eta}_\pm(v)) ,dS(V_\pm(v))\rangle
\tilde{\eta}_\pm(v)\,dv,
\end{multline*}
with $\tilde{\eta}_\pm(t)=e^{J \int_0^t
\left(E(U(v))-E(U_{\pm\infty})\right)\;dv}\left(z(t)+g(U(t),z(t))\right)$, 
but using the previous time decay estimates. 

There are two differences: 

One is in the estimate of 
the $H^s_{-\sigma}$. In fact before using the time decay estimates 
for $e^{-\rmi t H}P_c(H)$, we split the space associated with the 
continuous spectrum in two parts : one associated with energy closed 
to the thresholds and one associated to the rest of the spectrum. 
In the first part, we use the fact that $\sigma >  3/2$ to estimate 
the $H^s_{-\sigma}$  by the $B^\beta_{\infty,2}$ norm since we work with 
bounded energies. In the second part, since we work far from thresholds, 
we use Proposition \ref{Prop:PropagationHighEnergy} after estimating 
the $H^s_{-\sigma}$ by the $H^s_{-3/2}$. 

The other difference is in the estimation of the $B^\beta_{\infty,2}$ norm. 
We use Corollary \ref{Cor:LinDispersion} for $e^{Jt
H(V_{\pm\infty})}z_0$ and Theorem \ref{Thm:Dispersion} for the integrals.
\end{proof}

We have that $\lim_{\pm\infty}U=U_{\pm \infty}$ exist. If $z_0\in H^s_\sigma$ then the
associated solution $U$ satisfies
\begin{equation*}
\left|\dot{U}\right|\leq \frac{C}{\langle t\rangle^3}\left\|z_0\right\|_{H^s_\sigma}
\end{equation*}
and we have
\begin{equation*}
\int_0^t \left(E(U(v))-E(U_{\pm\infty}\right)\;dv \to E_{\pm\infty}
\mbox{ as } t\to \pm \infty
\end{equation*}
for some real $E_{\pm\infty}$. We introduce
\begin{equation*}
V_\pm(t)=e^{-\rmi \int_0^t
\left(E(U(v))-E(U_{\pm\infty})\right)\;dv}U(t),
\end{equation*}
they have a limit as $t\to \pm \infty$ respectively as being
\begin{equation*}
    V_{\pm \infty}=e^{-\rmi E_{\pm\infty}}U_{\pm \infty}.
\end{equation*}
Then we notice that we can also obtain an asymptotic profile for
$e^{\rmi t H +\rmi t E(U_\infty)}z(t)$ if $z_0$ is localized. But we
prefer to obtain a scattering result with respect to $e^{J t
H(V_\infty))}$. We have the
\begin{lemma}
                                                                \label{Lem:OnNonLinearScattering}
If Assumptions \ref{assumption:1}--\ref{assumption:Resonance} hold.
Then there exist~$\e_0>0$, ~$\delta_0>0$ such that for any
$\e\in(0,\e_0)$, $\delta\in(0,\delta_0)$, ~$U_0\in B_{\C^2}(0,\e)$
and $z_0\in B_{H^{s}_\sigma}(0,\delta)\cap \mathcal{H}_c(U_0)$ and
for the solution~$z$ of \eqref{Eq:ForZ} (with $U$ the solution of
\eqref{Eq:forU}) given in Lemma \ref{Lem:LocalWellPosednessForZ} the
following limit
\begin{equation*}
z_{\pm\infty}=\lim_{t\to\pm\infty}e^{-J t H(V_{\pm\infty})}e^{J
\int_0^t \left(E(U(v))-E(U_{\pm\infty})\right)\;dv}z(t)
\end{equation*}
exists in~$H^{s}$. Moreover, we have $z_{\pm\infty}\in {\mathcal
H}_c(V_{\pm\infty})\cap H^s_{\sigma}$ and there exists $C>0$ such
that
\begin{multline*}
\max\Big[\|e^{-J \int_0^t
\left(E(U(v))-E(U_{\pm\infty})\right)\;dv}e^{J t
H(V_{\pm\infty})\;dr}z_{\pm\infty}-z(t)\|_{H^{s}},\;\\ \|e^{-J \int_0^t
\left(E(U(v))-E(U_{\pm\infty})\right)\;dv}e^{J t
H(V_{\pm\infty})\;dr}z_{\pm\infty}-z(t)\|_{H^{s}_{-\sigma}},\;\\
\|e^{-J \int_0^t
\left(E(U(v))-E(U_{\pm\infty})\right)\;dv}e^{J t
H(V_{\pm\infty})\;dr}z_{\pm\infty}-z(t)\|_{B^\beta_{\infty,2}}
\Big]\leq
\frac{C}{\langle t \rangle^2}\left\|z_0\right\|_{H^{s}_\sigma}^2
\end{multline*}
and
\begin{equation*}
\|z_{\pm\infty}-e^{-J t
H(V_{\pm\infty})}e^{J \int_0^t
\left(E(U(v))-E(U_{\pm\infty})\right)\;dv}z(t)\|_{H^{s}_{3/2}}\leq \frac{C}{\langle t \rangle^{\frac{1}{2}}}\left\|z_0\right\|_{H^{s}_\sigma}^2.
\end{equation*}

\end{lemma}
\begin{proof}
Let
\begin{equation*}
t\mapsto \xi_\pm(t)=e^{J \int_0^t
\left(E(U(v))-E(U_{\pm\infty})\right)\;dv}z(t)
\end{equation*}
and
\begin{equation*}
t\mapsto V_\pm(t)=e^{-\rmi \int_0^t
\left(E(U(v))-E(U_{\pm\infty})\right)\;dv}U(t).
\end{equation*}
Using exactly the same method as the one of Lemma
\ref{Lem:StabilizationForT}, applied to
\begin{multline*}
e^{-Jt H(V_{\pm\infty})}\xi_\pm(t)=z_0 +\int_0^te^{-Js
H(V_{\pm\infty})}{\mathbf P}_c(V_\pm(v))J\left(d^2 F(S(V_\pm(v)))
-d^2 F(S(V_{\pm\infty}))\right)\xi_\pm(v)\,dv\\
+\int_0^te^{-Js H(V_{\pm\infty})}{\mathbf P}_c(V_\pm(v))JN(V_\pm(v),
\tilde{\eta}_\pm(v))\,dv\\
+\int_0^te^{-Js H(V_{\pm\infty})}{\mathbf
P}_c(V_\pm(v))dS(V(v))A(V_\pm(v),\tilde{\eta}_\pm(v))
\langle N(V_\pm(v),\tilde{\eta}_\pm(v)) ,dS(V_\pm(v))\rangle\,dv\\
-\int_0^te^{-Js H(V_{\pm\infty})}(d {\mathbf
P}_c(V_\pm(v)))A(V_\pm(v),\tilde{\eta}_\pm(v))\langle
N(V_\pm(v),\tilde{\eta}_\pm(v)) ,dS(V_\pm(v))\rangle
\tilde{\eta}_\pm(v)\,dv,
\end{multline*}
with $\tilde{\eta}_\pm(t)=e^{J \int_0^t
\left(E(U(v))-E(U_{\pm\infty})\right)\;dv}\left(z(t)+g(U(t),z(t))\right)$,
we prove that the limits
\begin{equation*}
    \lim_{t\to\pm\infty}e^{-Jt H(V_{\pm\infty})}\xi_\pm(t)=z_{\pm\infty}
\end{equation*}
exist. If we use the method of Lemma \ref{Lem:StabilizationForTloc},
we obtain the estimates on the convergence of $e^{Jt
H(V_{\pm\infty})}z_{\pm\infty}-\xi_\pm(t)$.
Then for multiplying by $e^{-J \int_0^t
\left(E(U(v))-E(U_{\pm\infty})\right)\;dv}$, we obtain  the
estimates and the convergence of
\begin{equation*}
    e^{-J \int_0^t
\left(E(U(v))-E(U_{\pm\infty})\right)\;dv}e^{-Jt
H(V_{\pm\infty})}z_{\pm\infty}-z(t).
\end{equation*}
Then since $\left(1-P_c(U(t))\right)z(t)=0$, we have
$\left(1-P_c(V_{\pm\infty})\right)z_{\pm\infty}=0$ and hence
$z_{\pm\infty}$ belongs to ${\mathcal H}_c(V_{\pm\infty})$.
\end{proof}
\section{The dynamic outside the center manifold}
                                                                \label{Sec:Outside}
We can make the same study in the center stable manifold and the
center unstable manifold but only in one direction of time. Let us
explain it for the center stable  manifold  in positive time since
it is similar for the center unstable manifold. Actually it is
equivalent if we revert the time direction.

We just give a sketch of the proof since it is similar to the
previous study. Using the idea of the proof of exponential
stabilization for Proposition \ref{Prop:CentrStabMan}, we write any
solution $\psi$ in the form $\phi + \rho + f(\phi,\rho)$ with $\phi$ in the center
manifold, $\rho\in X_s(V_0)$ and $f$ a function to be precised 
but ensuring that we are in the center stable manifold.

Indeed $W^{c}(V_0)$ is the graph of a smooth function
$h^c:X_{c}(V_0)\mapsto X_s(V_0)\oplus X_u(V_0)$ and $W^{cs}(V_0)$ 
is the graph of a smooth function
$h^u:X_{c}(V_0)\oplus X_s(V_0)\mapsto X_u(V_0)$. Let $\nu$ be such
that $\psi=S(V_0)+\nu$ satisfy \eqref{Eq:NLD}, we have
\begin{equation*}
\p_t \nu=JH(V_0)\nu+JN(V_0,\nu).
\end{equation*}
\begin{eqnarray*}
    \nu&=&y+h_c(y)+\rho+h^u(y,h^c(y)+\rho)\\
    &=&\phi(y)-S(V_0)+(\rho-\pi^s(V_0)h^c(y))+(h^u(y,h^c(y)+\rho)-\pi^u(V_0)h^c(y))\\
    &=&\phi(y)-S(V_0)+\rho+f(y,\rho)
\end{eqnarray*}
with $y=\pi^{c}(V_0)\nu=\pi^{c}(V_0)\left(\psi-S(V_0)\right)$ and $\phi(y)=S(V_0)+y+h^c(y)$ 
is in the center manifold 
and $\rho\in X_s(V_0)$. We have the following equation for $\rho$
\begin{equation}                                                \label{Eq:ForRho}
\p_t \rho= JH(V_0)\rho + M(V_0,y,\rho)
\end{equation}
where
\begin{multline*}
M(V_0,y,\rho)=\pi^{s}(V_0)\left\{JN(V_0,y+h^c(y)+\rho+f(y,\rho))-JN(V_0,y+h^c(y))\right\}\\
-\pi^{s}(V_0)Dh^c(y)\pi^{c}(V_0)\left\{JN(V_0,y+h^c(y)+\rho+f(y,\rho))-JN(V_0,y+h^c(y))\right\}.
\end{multline*}
Then we obtain for $\phi$ the equation
\begin{equation*}
    \p_t \phi= JH\phi + J\nabla F(\phi) + R(\phi,\rho)
\end{equation*}
\begin{equation*}
    R(\phi,\rho)=J\nabla F(\phi+\rho+f(y,\rho))-J\nabla F(\phi)-Jd^2 F(S(V_0))\rho - M(V_0,\pi^{c}(V_0)(\phi-S(V_0)),\rho)
\end{equation*}
with notice that $\left|R(\phi,\rho)\right|\leq
C(\left\|\phi\right\|_{H^s},\left\|\rho\right\|_{L^\infty})|\rho|$.

Working like in \ref{Sec:Stabilization}, we write $\phi= S(U)+\eta$ with 
$\eta=z+g\left(U,z\right)$ and we have the following equations for
$U$ and $z$:
\begin{equation*}
\begin{cases}
\dot{U}=-A(U,\eta)\langle N(U,\eta)-JR(U,\eta,\rho) ,dS(U)\rangle \\
\p_t z=JH(U)z+{\mathbf P}_c(U)JN(U,\eta)+{\mathbf
P}_c(U(v))dS(U(v))A(U(v),\eta(v)) \langle N(U(v),\eta(v))\\
\qquad -JR(U,\eta,\rho),dS(U(v))\rangle+(dP_c(U))A(U,\eta)\langle
N(U,\eta)-JR(U,\eta,\rho) ,dS(U)\rangle\eta+{\mathbf
P}_c(U)R(U,\eta,\rho)
\end{cases}
\end{equation*}
with
\begin{equation*}
\eta(t)=z(t) +g\left(U(t),z(t)\right).
\end{equation*}
where $g$ is defined by Lemma
\ref{Lem:DecompositionCentralManifofd2} and
\begin{equation*}
R(U,\eta,\rho)=R(S(U)+\eta,\rho)
\end{equation*}

These equations are similar to those we have studied but with an
extra term coming from $R$ which is exponentially decaying in
positive time. Indeed, let us introduce for any $T_0<0$ and $\gamma
\in(0,\gamma(V_0))$ and  $\delta>0$ the set
\begin{equation*}
{\mathcal R}_{T_0,\gamma}(\delta)=\left\{\rho\in{\mathcal
C}((T_0,+\infty),X_s(V_0)),\; \left|\rho(t)\right|_{H^{s}}\leq
\delta e^{-\gamma t}, \forall t>T_0\right\},
\end{equation*}
we study Equation \eqref{Eq:ForRho} in ${\mathcal
R}_{T_0,\gamma}(\delta)$ with small initial condition $\rho_0$.
We also define for any~$\e>0$
\begin{equation*}
{\mathcal U}_{T_0}(\e,\delta)\\
=\left\{U\in{\mathcal
C}^{1}((T_0,+\infty),B_{\C^2}(V_0,\e)),\;\|\dot{U}\|_{L^1((T_0,+\infty))\cap
L^\infty((T_0,+\infty))}\leq\delta^2\right\}
\end{equation*}
and for any~$U\in {\mathcal U}_{T_0}(\e)$, let~$s,\beta$ be such
that~\underline{$s> \beta+2>2$ and~$\sigma>3/2$},
\begin{multline*}
{\mathcal Z}_{T_0}(U,\delta)=\Bigg\{z\in {\mathcal C}((T_0,+\infty),
L^2(\R^3,\R^8)),\;z(t)\in {\mathcal
H}_c(U(t)),\,\\\max\left[\|z\|_{L^\infty((T_0,+\infty),H^{s})},
\|z\|_{L^2((T_0,+\infty),H^{s}_{-\sigma})},
\|z\|_{L^2((T_0,+\infty),B^\beta_{\infty,2})}\right]\leq\delta\Bigg\},
\end{multline*}
and $\e,\delta$ are small enough to ensure that for $U\in{\mathcal
U}(\e,\delta)$ and $z\in {\mathcal Z}(U,\delta)$
\begin{equation*}
    \ds S(U)+z+g(U,z)\in W^c(V_0)\cap B_{H^{s}}(S(V_0),r(V_0)).
\end{equation*}

For a sufficiently small $T_0$, we solve the equation for $z$ first
and then the one for $\rho$ and eventually the one for $U$ using the
method of Section \ref{Sec:Stabilization}. This gives us the desired
exponential decay for $\rho$ as well as similar results for $U$ and
$z$.

We notice that instead of Lemma \ref{Lem:LipshitzPropertyOfT}, we
obtain the
\begin{lemma}
If Assumptions \ref{assumption:1}--\ref{assumption:Resonance} hold.
Then for any $T>0$, there exist $T_0>0$, ~$\e_0>0$, ~$\delta_0>0$,
$C>0$ and $\kappa\in(0,1)$ such that for any $\e\in(0,\e_0)$,
$\delta\in(0,\delta_0)$, $U,\,U'\in {\mathcal U}_{T_0}(\e,\delta)$,
$\rho,\rho'\in{\mathcal R}_{T_0,\gamma}$ $z_0\in {\mathcal
H}_c(U(0))$, $z_0'\in {\mathcal H}_c(U'(0))$, $z\in {\mathcal
Z}_{T_0}(U,\delta)$ and $z'\in {\mathcal Z}_{T_0}(U',\delta)$, one
has
\begin{multline*}
\left\|z[z_0',U',\rho']
-z[z_0,U,\rho]\right\|_{L^\infty((T_0,T),H^{s})\cap L^2((T_0,T),L^{\infty})\cap L^2((T_0,T),H^{s}_{-\sigma})}\\
\leq C \left\|z_0-z_0'\right\|_{H^{s}} +\kappa\Big\{
\left\|U-U'\right\|_{L^\infty((T_0,T),\C^2}+\left\|\dot{U}-\dot{U}'\right\|_{L^\infty((T_0,T),\C^2)}\\
+\left\|e^{\gamma
t}(\rho-\rho')(t)\right\|_{L_t^\infty((T_0,T),X^s(V_0))}\Big\}.
\end{multline*}
\end{lemma}
Then for $\rho$ as a function of $U$, $z_0$ and $\rho_0$ (the
initial condition for $\rho$), we obtain the
\begin{lemma}
If Assumptions \ref{assumption:1}--\ref{assumption:Resonance} hold.
Then for any $T>0$ there exist $T_0>0$,~$\e_0>0$, ~$\delta_0>0$,
$C>0$ and $\kappa\in(0,1)$ such that for any $\e\in(0,\e_0)$,
$\delta\in(0,\delta_0)$, $U,\,U'\in {\mathcal U}_{T_0}(\e,\delta)$,
$r_0,r_0'\in X_s(V_0)$ $z_0\in {\mathcal H}_c(U(0))$, $z_0'\in
{\mathcal H}_c(U'(0))$, $z\in {\mathcal Z}_{T_0}(U,\delta)$ and
$z'\in {\mathcal Z}_{T_0}(U',\delta)$, one has
\begin{multline*}
\left\|e^{\gamma t}\left(\rho[z_0',U',\rho_0']
-\rho[z_0,U,\rho_0]\right)\right\|_{L_t^\infty((T_0,T),X^s(V_0))}\\
\leq C \left\|z_0-z_0'\right\|_{H^{s}} +\kappa\left\{
\left\|U-U'\right\|_{L^\infty((T_0,T),\C^2}+\left\|\dot{U}-\dot{U}'\right\|_{L^\infty((T_0,T),\C^2)}
+\left\|\rho_0-\rho_0'\right\|_{L^2}\right\}.
\end{multline*}
\end{lemma}

We also notice that the proof gives the wellposedness of Equation
\eqref{Eq:ForRho} in ${\mathcal R}_{T_0,\gamma}(\delta)$ with small
initial condition $\rho_0$ and that there exists $C>0$ such that the
solution $\rho$ satisfies
\begin{equation*}
    \left\|\rho(t)\right\|_{H^{s}}\leq C \left\|
\rho_\pm(0)\right\|e^{-\gamma t}, \forall t>T_0.
\end{equation*}
The asymptotic behaviour of $U$ and $z$ are obtained like in the previous section when $z_0$ is localized.

\section{End of the proof of main theorems}
                                                                    \label{Sec:EndProof}
We notice that the small {\it locally invariant} center manifold
build in Section \ref{Sec:CSUMan} for Equation \eqref{Eq:Centered}
is now a small {\it invariant} (globally in time) center manifold.
Indeed, we have just proved the stabilization towards the PLS
manifold, this ensures that a solution in the center manifold will
stay inside this manifold in the two direction of time.

Now let us consider ${\mathcal CM}$ as being the union of all
these small globally invariant center manifolds and $0$. Using the uniqueness
of center manifold and Lemma
\ref{Lem:DecompositionCentralManifofd2}, we prove that ${\mathcal CM}\setminus \{0\}$ is a
manifold. Now we generalize Lemma
\ref{Lem:DecompositionCentralManifofd2} by the
\begin{lemma}																												\label{Lem:DecompositionCentralManifofd3}
For any $s,s',\sigma\in\R$ and $p,q\in[1,\infty]$, there exist
$\e>0$, a continuous map $r:B_\C^2(0,\e)\mapsto \R^+$ with
$r(U)=O(\Gamma(U))$ and a continuous map $\Psi: S\mapsto {\mathcal
CM}$ where
\begin{equation*}
    {\mathcal S}_\sigma=\left\{(U,z);\;U\in B_{\C^2}(0,\e),z\in{\mathcal H}_c(U)\cap
B_{H^{s'}_\sigma}(0,r(U))\right\}
\end{equation*}
is endowed with the metric of $\C^2\times H^{s'}_\sigma$.

Moreover $\Psi$ is bijective from ${\mathcal S}$ to an open
neigborhood of $(0,0)$ in ${\mathcal CM}$ and smooth on ${\mathcal
S}\setminus\{(0,0)\}$. For all $U\in B_{\C^2}(0,\e)$, there exists
$C>0$ such that for all $z\in{\mathcal H}_c(U)\cap
B_{H^{s'}}(0,r(U))$, $\Psi(U,z)\in {\mathcal H}_1(U)$,
$z+\Psi(U,z)\in{\mathcal H}_0(U)^\bot$, $S(U)+z+\Psi(U,z)\in
{\mathcal CM}$. For sufficiently small non zero $U$, we have
$\left\|\Psi(U,z)\right\|_{B^s_{p,q}}=O(
\left\|z\right\|_{H^{s'}}^2)$ for $z\in H^{s'}$ such that $(U,z)\in
{\mathcal S}$.
\end{lemma}
\begin{proof}
The proof works like for Lemma
\ref{Lem:DecompositionCentralManifofd2}. 
The statements for $r$ follow from Remark \ref{Rem:TauGamma}.

\end{proof}

The scattering result follows from a one to one correspondence of
the initial profile with the asymptotic profile as stated in the
\begin{proposition}
If Assumptions \ref{assumption:1}--\ref{assumption:Resonance} hold.
There exist $\e>0$  and a continuous map $r:B_\C^2(0,\e)\mapsto
\R^+$ with $r(U)=O(\Gamma(U))$ and ${\mathcal V}_\sigma$,  ${\mathcal
V}_\pm$ neighborhoods of $(0,0)$ in
\begin{equation*}
    S_\sigma=\left\{\left(U,z\right);\;U\in \C^2,\; z\in {\mathcal H}_c(U)\cap B_{H^{s}_\sigma}(0,r(U)) \right\}
\end{equation*}
endowed with the norm of $\C^2\times H^{s}_\sigma$ such that the
maps
\begin{equation*}
   \mathcal{P}_\pm:
   \left(\begin{array}{c}U_0\\z_0\end{array}\right) \in {\mathcal V}_\sigma
   \mapsto  \left(\begin{array}{c}V_{\pm\infty}\\z_{\pm\infty}\end{array}\right) \in {\mathcal V}_\pm \\
\end{equation*}
are bijections and are smooth on ${\mathcal V}_0\setminus\{(0,0)\}$.
\end{proposition}
\begin{proof}
We choose for example
\begin{equation*}
    {\mathcal V}_\sigma=\left\{\left(U,z\right);\;U\in B_{\C^2}(0,\e),\; z\in {\mathcal H}_c(U)\cap B_{H^{s}_\sigma}(0,r(U)) \right\}
\end{equation*}
for some positive $\e$ and we work on the manifold ${\mathcal
V}_\sigma\setminus\left\{\right(0,0)\}$ which is locally isomorphic to an
open set of $\C^2\times {\mathcal H}_c(U)\cap H^{s}_\sigma$. We write
\begin{equation*}
\mathcal{P}_\pm^{U_0}\left(U,z\right)=\left(U,z\right)+\mathcal{R}_\pm^{U_0}\left(U,z\right)
\end{equation*}
Since
\begin{equation*}
\left\|\left(U_\infty,z_\infty\right)-\left(U_0,z_0\right)\right\|_{H^{s}_\sigma}=O\left(|U_0|^2+\left\|z_0\right\|_{H^{s}_\sigma}^2\right),
\end{equation*}
we only need to prove the statement locally. Hence we prove that in
a neighborhood of $(U_0,0)$. The maps
$\mathcal{P}_\pm^{U_0}\left(U,z\right)\mapsto(Id_{\C^2},P_c(U_0))\mathcal{P}_\pm\left(U,R(U,U_0)z\right)$
are bijective ($P_c$ and $R$ are defined in Proposition
\ref{Prop:ContSpectrumLinOp}).

To prove that $\mathcal{P}_\pm^{U_0}$ is bijective ({\it i.e.} the scattering exists). Let us prove it for 
$\mathcal{P}_+^{U_0}$ (it is similar fo $\mathcal{P}_-^{U_0}$). It is enough to prove that the following system has a unique solution in 
an open neighborhood of $(0,0)$ in ${\mathcal S}_\sigma$: 
\begin{equation*}
V_\pm(t)=V_{\pm\infty}+\int_t^\infty A(V_\pm(v),e^{Js
H(V_{\pm\infty})}\tilde{\eta}_\pm(v))\langle N(U(v),e^{Js
H(V_{\pm\infty})}\tilde{\eta}_\pm(v))
,dS(V_\pm(v))\rangle\,dv,
\end{equation*}
and
\begin{multline*}
\tilde{\xi}_+(t)=z_\infty -\int_t^\infty e^{-Js
H(V_{+\infty})}{\mathbf P}_c(V_+(v))J\left(d^2 F(S(V_+(v)))
-d^2 F(S(V_{+\infty}))\right)e^{Js
H(V_{+\infty})}\tilde{\xi}_+(v)\,dv\\
-\int_t^\infty e^{-Js H(V_{+\infty})}{\mathbf P}_c(V_+(v))JN(V_+(v),
e^{Js
H(V_{+\infty})}\tilde{\eta}_+(v))\,dv\\
-\int_t^\infty e^{-Js H(V_{+\infty})}{\mathbf
P}_c(V_+(v))dS(V(v))A(V_+(v),e^{Js
H(V_{+\infty})}\tilde{\eta}_+(v))
\langle N(V_+(v),e^{Js
H(V_{+\infty})}\tilde{\eta}_+(v)) ,dS(V_+(v))\rangle\,dv\\
+\int_t^\infty e^{-Js H(V_{+\infty})}(d {\mathbf
P}_c(V_+(v)))A(V_+(v),e^{Js
H(V_{+\infty})}\tilde{\eta}_+(v))\langle
N(V_+(v),e^{Js
H(V_{+\infty})}\tilde{\eta}_+(v)) ,dS(V_+(v))\rangle
e^{Js
H(V_{+\infty})}\tilde{\eta}_+(v)\,dv,
\end{multline*}
with $\tilde{\eta}_+(t)=\tilde{\xi}_+(t)+e^{-Js
H(V_{+\infty})}g\left(V_{+}(t),e^{Js
H(V_{+\infty})}\tilde{\xi}_+(t))\right)$.

This system can be solved by a fixed point argument in the set of function such that
\begin{multline*}
\max\Big[\sup_{t\in \R}(\|z_{+\infty}-\tilde{\xi}_+(t)\|_{H^{s}}),\; \sup_{t\in
\R}\langle t\rangle^{3/2}\|z_{+\infty}-\tilde{\xi}_+(t)\|_{H^{s}_{-\sigma}},\;\\
\sup_{t\in \R}\langle t\rangle^{3/2}\|z_{+\infty}-\tilde{\xi}_+(t)\|_{B^\beta_{\infty,2}},\;
\sup_{t\in \R}(\langle t\rangle^{-3/2}\|z_{+\infty}-\tilde{\xi}_+(t)\|_{H^{s}_{3/2}})\Big]
\end{multline*}
and
\begin{equation*}
\left\langle t \right\rangle^2\left|V_+(t)-V_{+\infty}\right|	
\end{equation*}
are small with the method we used in Lemma \ref{Lem:OnNonLinearScattering}.


\end{proof}

For the same reasons the small locally invariant center-stable
manifold build in Section \ref{Sec:CSUMan} is invariant in positive
time. We can also consider the union of these manifolds, and we can
obtain a map $\Phi_+$ similar to the map $\Psi$ built in Lemma 
\ref{Lem:DecompositionCentralManifofd3}.
The instability in negative time is in fact a consequence of
Proposition \ref{Prop:CentrUnstabMan}.

The corresponding conclusion holds for the center unstable manifold.

The statements on the instability outside these manifolds follow
from Propositions \ref{Prop:CentrStabMan} and
\ref{Prop:CentrUnstabMan}.

\section*{\bf Acknowledgements}
I would like to thank \'Eric S\'er\'e for fruitful discussions and
advices during the preparation of this work. I am also indebted to
Galina Perelman for her careful reading, her comments and
suggestions.

\bibliographystyle{alpha}
\bibliography{biblio}

\end{document}